\documentclass[12pt]{article}
\usepackage{amsmath, amsfonts, amssymb,amsthm, amscd}

\usepackage{graphicx}
\usepackage{psfrag}
\usepackage{color}
\usepackage{hyperref}

\theoremstyle{definition}

\theoremstyle{remark}

\newtheoremstyle{thm}
  {12pt}
  {12pt}
  {\itshape}
  {\parindent}
  {\scshape}
  {.}
  {5pt}
  {}

\theoremstyle{thm}
\newtheorem*{T6.17}{Theorem \ref{gradientthm1}}
\newtheorem{thm}{Theorem}[section]

\newtheoremstyle{prop}
  {12pt}
  {12pt}
  {\itshape}
  {\parindent}
  {\scshape}
  {.}
  {5pt}
  {}

\theoremstyle{prop}
\newtheorem{prop}[thm]{Proposition}

\newtheoremstyle{lem}
  {12pt}
  {12pt}
  {\itshape}
  {\parindent}
  {\scshape}
  {.}
  {5pt}
  {}

\theoremstyle{lem}
\newtheorem*{L2.6}{Lemma \ref{lem2.4}}
\newtheorem{lem}[thm]{Lemma}

\newtheoremstyle{defn}
  {12pt}
  {12pt}
  {\itshape}
  {\parindent}
  {\scshape}
  {.}
  {5pt}
  {}

\theoremstyle{defn}
\newtheorem{defn}[thm]{Definition}

\newtheoremstyle{examp}
  {12pt}
  {12pt}
  {}
   {\parindent}
  {\scshape}
  {.}
  {5pt}
  {}

\theoremstyle{examp}

\newtheoremstyle{cor}
  {12pt}
  {12pt}
  {\itshape}
  {\parindent}
  {\scshape}
  {.}
  {5pt}
  {}

\theoremstyle{cor}

\newtheoremstyle{recipe}
  {12pt}
  {12pt}
  {\itshape}
   {\parindent}
  {\scshape}
  {.}
  {5pt}
  {}

\theoremstyle{recipe}

\newtheoremstyle{rem}
  {12pt}
  {12pt}
  {}
   {\parindent}
  {\scshape}
  {.}
  {5pt}
  {}

\theoremstyle{rem}
\newtheorem{rem}[thm]{Remark}

\newcommand{\bp}{\begin{proof}}
\newcommand{\ep}{\end{proof}}

\newcommand{\norm}[1]{\left\Vert#1\right\Vert}
\newcommand{\abs}[1]{\left\vert#1\right\vert}

\newcommand{\ssc}{\text{sc}}

\renewcommand{\epsilon}{\varepsilon}

\newcommand{\what}{\widehat}
\newcommand{\wh}{\widehat}

\newcommand{\ov}{\overline}

\newcommand{\ind}{\operatorname{Ind}}

\newcommand{\R}{{\mathbb R}}

\providecommand{\ker}[1]{$\text{ker}\ {#1}$}
\newcommand{\trl}{\triangleleft}
\newcommand{\N}{{\mathbb N}}

\def\abs#1{\mathopen|#1\mathclose|}

\def\norm#1{\mathopen\|#1\mathclose\|}

{\catcode`"=12 \gdef\hex{"}}

\mathchardef\laplace=\hex0001

\mathchardef\nabla=\hex0272

\makeatletter

\def\@@dalembert#1#2{\setbox0\hbox{$#1\mathrm I$}

  \vrule height\ht0 depth\z@ width.04\ht0

  \rlap{\vrule height\ht0 depth-.96\ht0 width.8\ht0}

  \vrule height.1\ht0 depth\z@ width.8\ht0

  \vrule height\ht0 depth\z@ width.1\ht0 }

\def\dalembert{\mathbin{\mathpalette\@@dalembert{}}\,}

\makeatother

\begin{document}

\title{A General Fredholm Theory I:\\
A Splicing-Based Differential Geometry}

\author{ H. Hofer\footnote{Research partially supported
by NSF grant  DMS-0603957.}\\
New York University\\ USA\and K. Wysocki \footnote{Research
partially supported by  NSF grant DMS-0606588. }\\Penn State University\\
USA\and E. Zehnder
 \footnote{Research partially supported by  TH-project.}\\ETH-Zurich\\Switzerland}

\maketitle

\section{Introduction}
In a series of papers we develop a generalized Fredholm theory and
demonstrate its applicability to a variety of problems including
Floer theory, Gromov-Witten theory, contact homology, and symplectic
field theory. Here are some of the basic common features:\\

\noindent $\bullet$  The moduli spaces are solutions of elliptic
PDE's showing serious non-compactness phenomena having well-known
names like bubbling-off, stretching the neck, blow-up, breaking of
trajectories. These drastic names  are a manifestation of the fact
that  one is confronted with analytical limiting phenomena where
the classical analytical descriptions break down.\\

\noindent $\bullet$ When the moduli spaces are not compact, they
admit nontrivial compactifications  like the Gromov
compactification, \cite{G}, of the space of pseudoholomorphic curves
in Gromov-Witten theory or the compactification of the moduli spaces
in symplectic field theory (SFT) as described in \cite{BEHWZ}.\\

\noindent $\bullet$   In many problems like in Floer theory, contact
homology or symplectic field theory the algebraic structures of
interest are precisely those created by  the ``violent analytical
behavior" and its  ``taming" by suitable compactifications. In fact,
the algebra is created by the complicated interactions of many
different moduli spaces.\\

In the abstract theory we shall introduce a new class of spaces
called polyfolds which in applications are the ambient spaces of the
compactified moduli spaces. We introduce bundles $p:Y\rightarrow X$
over polyfolds which, as well as the underlying polyfolds, can have
varying dimensions. We define the notion of a Fredholm section
$\eta$ of the bundles $p$ whose zero sets $\eta^{-1}(0)\subset X$
are in our applications precisely the compactfied moduli spaces one
is interested in. The normal ``Fredholm package'' will be
constructed consisting
 of an abstract perturbation and transversality theory.
 In the case of transversality the solution spaces are smooth manifolds,
 smooth orbifolds, or smooth weighted branched manifolds
 (in the sense of McDuff, \cite{Mc}),
 depending on the generality of the situation.

 The usefulness of this  theory will be
 illustrated by our `Application Series'.
The applications include Gromov-Witten theory, Floer theory and SFT,
see \cite{HWZ4,HWZ5}.
 It is, however,
clear that the theory applies to many more nonlinear problems
showing a lack of compactness.

The current paper is the first in the `Theory Series' and deals with
a generalization of differential geometry which is based on new
local models. These local models are open sets in splicing cores.
Splicing cores are smooth spaces with tangent spaces having  in
general locally varying dimension. These spaces are associated to
splicings which is the basic concept in this paper. The so obtained
local models for a new kind of smooth spaces  are needed to deal
with the functional analytical descriptions of situations in which
serious compactness problems arise. We would also like to note that
the applications of the concepts in this paper can be viewed as a
generalization of \cite{El} to a situation where we have varying
domains and targets.

The second paper, \cite{HWZ2}, develops the implicit function
theorems in this general context and extends the usual Fredholm
theory.

The third paper, \cite{HWZ3}, develops the Fredholm theory in
polyfolds, which could be viewed as a theory of Fredholm functors in
a version of Lie groupoids with object and morphism spaces build on
the new local models (see \cite{Mj},\cite{MM} for the groupoid
concepts in a manifold world).  The Fredholm theory in this
generalization is sufficient to deal with the problems mentioned
above.

On purpose we have not included any applications in this series
since we did not want to dilute the ideas. The conceptual framework
 should apply to many more situations. We refer the reader to
\cite{HWZ-polyfolds1,HWZ-polyfolds2,HWZ4,HWZ5} for applications on
different depth levels. An overview is given in \cite{Hofer}.\\

\noindent {\bf Acknowledgement:} We would like to thank P. Albers
and U. Hryniewicz for helpful comments.

\section{Sc-Calculus in Banach Spaces}\label{chapter1} In
order to develop the generalized nonlinear Fredholm theory needed
for the symplectic field theory, we  start with calculus issues. In
a first step we equip  Banach spaces with the structure of a scale,
called  sc-structure. Scales are a well-known concept from
interpolation theory, see for example \cite{Tr}. We give a new
interpretation of a scale as a generalization of a smooth structure.
Then we introduce the appropriate class of smooth maps. Having
developed the notion of an sc-smooth structure on an open subset of
a Banach space as well as that of a smooth map, the validity of the
chain rule allows then, in principle, to develop an
``sc-differential geometry'' by simply imitating the classical
constructions. However, new objects are possible, with the most
important one being that of a general splicing. The main purpose of
this paper is to introduce them and to show how they define  local
models for a new class of smooth spaces, which are crucial for the
afore-mentioned applications.

\subsection{sc-Structures}\label{sect1.1} We begin by introducing
the notion of an  sc-smooth structure on a Banach space and on its
open subsets.
\begin{defn}\label{scstructure}
Let $E$ be a Banach space.  An {\bf  sc-structure} \index{sc-smooth
structure} on $E$ is given by a nested sequence
$$E=E_0 \supseteq E_1\supseteq E_2\supseteq \cdots \supseteq
\bigcap_{m\geq 0}E_m=:E_{\infty}$$ of Banach spaces $E_m$, $m\in
{\mathbb N}=\{0,1,2,\cdots\}$, having the following properties.
\begin{itemize}
\item[$\bullet$] \:  If $m<n$, the inclusion $E_n\hookrightarrow E_m$  is a
compact operator.

\item[$\bullet$] \:The vector space $E_{\infty}$
is dense in $E_m$  for every $m\geq 0$.
\end{itemize}
\end{defn}
In the following we shall sometimes talk about an sc-smooth
structure on a Banach spaces rather than an sc-structure to
emphasize the smoothness aspect.  From the definition of an
sc-structure it follows, in particular, that $E_n\subseteq E_m$ is
dense if $m<n$  and the embedding is continuous . We note that
$E_{\infty}$ has the structure of a Frechet space. In the case
$\dim(E)<\infty$ the only possible sc-structure is the constant
structure with $E_m=E$.

If $U\subset E$ is an open subset we define the {\bf induced
sc-smooth structure}\index{induced sc-structure} on $U$  to be the
nested sequence $U_m=U\cap E_m$. Given an sc-smooth structure on $U$
we observe that $U_m$ inherits the sc-smooth structure defined by
${(U_m)}_k=U_{m+k}$. We will write $E^m$ to emphasize that we are
dealing with the Banach space $E_m$ equipped with the sc-structure
$(E^m)_{k}:=E_{m+k}$ for all $k\geq 0$. Similarly we will
distinguish between  $U_m$ and $U^m$.
\begin{rem}\label{remake}
The compactness requirement is crucial for applications. It is
possible to develop a theory without this requirement, but it is not
applicable to the theories we are interested in. This alternative
theory would in the case of the constant sequence $E_m=E$ recover
the standard smooth structure on $E$. However, the notion of
 a smooth map would be more restrictive. Both theories, the one
 described in this paper and the one just alluded to, intersect therefore only
 in
 the standard finite-dimensional theory. See Remark \ref{axx} for
 further details.
\end{rem}

If $E$ and $F$ are equipped with sc-structures, the Banach space
$E\oplus F$ carries the  sc-structure defined by $(E\oplus
F)_m=E_m\oplus F_m$.

\begin{defn}\label{def1.3}
Let $U$ and $V$ be open subsets of sc-smooth Banach spaces. A
continuous map $\varphi:U\rightarrow V$ is said to be of
$\mathbf{class}\ \mathbf{sc}^{\mathbf{0}}$ \index{class! sc$^0$} or
simply  $ \mathbf{sc}^{\mathbf{0}}$  if $\varphi(U_m)\subset V_m$
and the induced  maps
$$
\varphi:U_m\rightarrow V_m
$$
are all continuous.
\end{defn}

Next we define the tangent bundle.
\begin{defn}
Let $U$ be an open subset in an  sc-smooth Banach space $E$ equipped
with the induced sc-structure. Then the {\bf tangent bundle} $TU$
of $U$ \index{tangent! bundle} is defined by $TU=U^1\oplus E$. Hence
the induced sc-smooth structure is defined  by  the nested sequence
$$(TU)_m=U_{m+1}\oplus E_m$$
 together with the sc$^0$-projection
$$
p:TU\rightarrow U^1.
$$
\end{defn}
Note that the tangent bundle is not  defined on $U$ but  merely on
the smaller, but dense, subset  $U_1$. We shall refer to the points
in $U_m$ sometimes as  points in $U$ on the level $m$.

\subsection{Linear Sc-Theory}\label{section2.2}
 We begin by developing some of the linear theory needed in the
 sc-calculus.
\begin{defn}
Consider $E$ equipped with an sc-smooth structure.
\begin{itemize}
\item An  {\bf sc-subspace}  $F$ of $E$ consists of a closed linear
subspace $F\subseteq E$, so that $F_m= F\cap E_m$ defines an
sc-structure for $F$.
\item An sc-subspace $F$ of $E$  {\bf splits}  if  there exists
another sc-subspace $G$ so that on every level we have the
topological direct sum
$$
E_{m}=F_{m}\oplus G_{m}.
$$
\end{itemize}
\end{defn}
We  shall use the notation $ E= F\oplus_{sc} G $
 or
 $
 E=F\oplus G
 $
 if there is no possibility of confusion.

Next we introduce the relevant linear operators in the sc-context.
\begin{defn}
Let $E$ and $F$ be sc-smooth Banach spaces.
\begin{itemize}
\item An  {\bf sc-operator} \index{sc-operator} $T:E\to F$
is a bounded linear operator which in addition is $\ssc^0$, i.e.~it
induces bounded linear operators $ T:E_{m}\rightarrow F_{m}$ on all
levels.

\item  An {\bf sc-isomorphism} \index{sc-isomorphism}
is a  bijective sc-operator $T:E\rightarrow F$  such that
$T^{-1}:F\rightarrow E$ is also an  sc-operator.
\end{itemize}
\end{defn}
An interesting class of sc-operators is the class of sc-projections,
i.e. sc-operators $P$ with $P^2=P$.

\begin{prop}\label{split} Let $E$ be an sc-smooth Banach space and
$K$ a finite-dimensional subspace of $E_{\infty}$. Then $K$ splits
the sc-space $E$.
\end{prop}
Note  that a finite-dimensional subspace $K$ of $E$ which splits the
sc-smooth space $E$  is necessarily  a subspace of $E_{\infty}$.
\begin{proof}
Take a basis $e_1,...,e_n$ for $K$ and fix the associated dual
basis. By Hahn Banach this dual basis can be extended to continuous
linear functionals $\lambda_1,...,\lambda_n$ on $E$. Now
$P(h)=\sum_{i=1}^{n} \lambda_i(h)e_i$ defines a continuous
projection on $E$ with image in $K\subset E_{\infty}$. Hence $P$
induces continuous maps $E_m\rightarrow E_m$. Therefore $P$ is an
sc-projection. Define $Y_m=(Id-P)(E_m)$.  Setting  $Y=Y_0$  we have
$E=K\oplus Y$.   By construction, $Y_m\subset E_m\cap Y_0$.  An
element  $x\in E_m\cap Y_0$ has the form $x=e-P(e)$ with  $e\in
E_0$. Since $P(e)\in E_{\infty}$ we see that $e\in E_m$,  implying
$Y_m = E_m\cap Y_0$. Finally, $Y_{\infty}=\cap_{m\geq 0}Y_m$ is
dense in $Y_m$ for every $m\geq 0$.  Indeed, if $x\in Y_m$, we can
choose $x_k\in E_{\infty}$ satisfying $x_k\rightarrow x$ in $E_m$.
Then $(Id-P)x_k\in Y_{\infty}$ and $(Id-P)x_k \rightarrow (Id-P)x =
x$ in $Y_m$.
\end{proof}

 We can introduce the notion of a linear Fredholm
operator  in the sc-setting.
\begin{defn}\label{def1.18}
Let $E$ and $F$ be sc-smooth Banach spaces. An sc-operator
$T:E\rightarrow F$ is called  {\bf Fredholm} \index{sc-Fredholm
operator} provided there exist sc-splittings $E=K\oplus_{sc} X$ and
$F=Y\oplus_{sc} C$ having the following properties.
\begin{itemize}
\item $K=\text{kernel\ $(T)$}$ is finite-dimensional.
\item $C$ is finite-dimensional.
\item  $Y=T(X)$  and $T:X\rightarrow Y$ defines a linear
sc-isomorphism.
\end{itemize}
\end{defn}

The above definition implies that $T(X_m)=Y_m$,  the kernel of
$T:E_m\rightarrow F_m$ is equal to $K$ and  $C$ spans its cokernel,
so that
$$
E_m=K \oplus X_m\:\;  \text{and}\: \; F_m=C\oplus  T(E_m)
$$
for all $m\geq 0$. It is an easily established fact that the
composition of two sc-Fredholm operators $T$ and $S$  is
sc-Fredholm. From the index additivity of classical Fredholm
operators we obtain the same in our set-up
$$
i(TS)=i(T)+i(S).
$$
The following observation, called the regularizing property, should
look familiar.
\begin{prop}\label{regular}
Assume $T:E\rightarrow F$ is sc-Fredholm and
$
T(e)\in F_m
$
for some  $e\in E_0$. Then $e\in E_m$.
\end{prop}
\begin{proof}
Since $F_m=T(E_m)\oplus C$, the element $f=T(e)\in F_m$ has the
representation
$$f=T(x)+c$$
for some $x\in X_m$ and $c\in C$. Similarly, $e$ has the
representation
$$e=k+x_0,$$
with  $k\in K=\ker T$ and $x_0\in X_0$ because $E_0=K\oplus X_0$.
From $T(e)=f=T(x)+c$ and $T(x)=T(x_0)$ one concludes $T(x_0-x)=c$.
Hence $c=0$  because $T(E_0)\cap C=\{0\}$.  Consequently,  $x_0-x\in
K$. Since $e-x=k+(x_0-x)\in K$ and $x\in E_m$ and $K\subset E_m$,
one  concludes $e\in E_m$ as claimed.
\end{proof}

We end this subsection with the important definition of an
$\ssc^+$-operator and an stability result for Fredholm maps.
\begin{defn}\label{def.1.20}
Let $E$ and $F$ be sc-Banach spaces. An sc-operator $R:E\rightarrow
F$ is said to be an
$\mathbf{sc}^{\boldsymbol+}${\bf-}$\mathbf{operator}$\index{sc$^+$-operator}
if  $R(E_m)\subset F_{m+1}$ for every $m\geq 0$
and if $R$  induces
an $\ssc^0$-operator $E\rightarrow F^1$.
\end{defn}
Let us note that due to the (level-wise) compact embedding
$F^1\rightarrow F$ an  sc$^+$-operator induces on every level a
compact operator. This follows immediately from the factorization
$$
R:E\rightarrow F^1\rightarrow F.
$$
The stability  result is  the following statement.
\begin{prop}\label{prop1.21}
Let $E$ and $F$ be sc-Banach spaces. If $T:E\rightarrow F$ is an
sc-Fredholm operator and $R:E\rightarrow F$ an sc$^+$-operator, then
$T+R$ is also an  sc-Fredholm operator.
\end{prop}
\begin{proof}
Since $R:E_m\rightarrow F_m$ is compact for every level we see that
$T+R:E_m\rightarrow F_m$ is Fredholm for every $m$. Let $K_m$ be the
kernel of $T+R:E_m\rightarrow F_m$. We claim that $K_{m}=K_{m+1}$
for every $m\geq 0$. Clearly, $K_{m+1}\subseteq K_{m}$. To see that
$K_m\subseteq K_{m+1}$, take  $x\in K_m$. Then  $Tx=-Rx\in F_{m+1}$
and, in view of  Proposition \ref{regular},   $x\in E_{m+1}$. Hence
$x\in K_{m+1}$ implying  $K_m\subseteq K_{m+1}$. Set $K=K_0$.  By
Proposition \ref{split}, $K$ splits the sc-space $E$  since it is a
finite dimensional subset of $E_{\infty}$. Hence we have the
sc-splitting $E=K\oplus X$  for a suitable sc-subspace  $X$. Next
define $Y_m=(T+R)(E_m)$. This defines an  sc-structure on $Y=Y_0$.
Let us show that $F$ induces an sc-structure on $Y$ and that this is
the one given by $Y_m$. For this it suffices  to show that
\begin{equation}\label{eqscp}
Y\cap F_m = Y_m.
\end{equation}
Clearly,
$$
Y_m = (T+R)(E_m)= F_m\cap (T+R)(E_m)\subset F_m\cap (T+R)(E_0) =
Y\cap F_m.
$$
Next  assume that $y\in Y\cap F_m$. Then there exists $x\in E_0$
with $Tx+Rx =y$. Since $R$ is an  sc$^+$-section it follows that
$y-Rx\in F_1$ implying that $x\in E_1$. Inductively we find that
$x\in E_m$ implying that $y\in Y_m$ and \eqref{eqscp} is proved.
Observe that we also have
$$
F_{\infty}\cap Y= (\bigcap_{m\in{\mathbb N}} F_m)\cap
Y=\bigcap_{m\in {\mathbb N}} (F_m\cap Y)= \bigcap_{m\in {\mathbb N}}
Y_m=Y_{\infty}
$$

In view of Lemma \ref{dense} below, there exists a finite
dimensional subspace $C\subset F_{\infty}$ satisfying $F_0=C\oplus
Y$.   From this it follows that $F_m=C\oplus Y_m$.  Indeed, since
$C\cap Y_m\subset C\cap Y$,  we have $C\cap Y_m=\{0\}$. If $f\in
F_m$, then $f=c+y$ for some $c\in C$ and $y\in Y$ since $F_m\subset
Y$ and $F_0=C\oplus Y$. Hence $y=f-c\in F_m$ and using  Proposition
\ref{regular} we conclude $y\in F_m$.  This implies, in view of
\eqref{eqscp}, that $F_m=C\oplus Y_m$.   We also have
$F_{\infty}=C\oplus (F_{\infty}\cap Y)=C\oplus Y_{\infty}$. It
remains to show that $Y_{\infty}$ is dense in $Y_m$ for every $m\geq
0$.  Take $y\in Y_m$. Then $y=(T+R)(x)$ for some $x\in E_m$. The
space  $E_{\infty}$ is dense in $E_m$   so that there  exists a
sequence $(x_n)\subset E_{\infty}$ converging to $x$ in $E_m$.  The
operator $T+R$ is $\ssc^0$-continuous and so the sequence
$y_n:=(T+R)(x_n)\in F_m$ converges to $y=(T+R)(x)$ in $F_m$.  Now,
in view of Proposition \ref{regular} and Definition \ref{def.1.20},
the points  $y_n$ belong to $Y_{\infty}$ and our claim is proved.
Consequently, we have the sc-splitting
$$
F= Y\oplus_{sc} C
$$
and, up to the Lemma \ref{dense} below, the proof of the proposition
is complete.
\end{proof}

\begin{lem}\label{dense}
Assume  $F$ is a  Banach space  and $F=D\oplus Y$ with $D$ of
finite-dimension and  $Y$ a closed subspace of $F$. Assume, in
addition, that $F_{\infty}$ is a dense subspace of $F$. Then there
exists a finite dimensional subspace $C\subset F_{\infty}$ such that
$F=C\oplus Y$.
\end{lem}
\begin{proof}
The quotient $F/Y$ is a finite-dimensional Banach space and we have
a continuous projection operator
$$
p:F\rightarrow F/Y
$$
Since $F_\infty$ is a dense linear subspace of $F$ we find that
$$
p(F_\infty)=F/Y.
$$
 Take any basis for $F/Y$ and pick representatives for these
vectors in $F_\infty$. Their span $C$ has the desired property.
\end{proof}

\subsection{Sc-Smooth Maps}\label{sec1.2}
In this subsection we introduce the notion of an  sc$^1$-map.
\begin{defn}\label{sscc1} Let $E$ and $F$
be sc-smooth Banach spaces and let $U\subset E$ be an open subset.
An  $\ssc^0$-map $f:U\rightarrow F$ is said to be
$\mathbf{sc}^{\mathbf{1}}$ or of $\mathbf{class}\
\mathbf{sc}^{\mathbf{1}}$  \index{class! sc$^1$} if  the following
conditions hold true.
\begin{itemize}\label{sc-1}
\item[(1)] For every $x\in U_1$ there exists a
linear map $Df(x)\in  {\mathcal L}(E_0, F_0)$ satisfying for $h\in
E_1$, with $x+h\in U_1$,
$$\frac{1}{\norm{h}_1}\norm{f(x+h)-f(x)-Df(x)h}_0\to 0\quad
\text{as\ $\norm{h}_1\to 0$}\mbox{}\\[4pt]$$
\item[(2)]  The {\bf tangent map} \index{tangent!map} $Tf:TU\to TF$,
defined by
$$Tf(x, h)=(f(x), Df(x)h)
$$
is an $\ssc^0$-map.
\end{itemize}
\end{defn}

The linear map $Df(x)$ will in the following often be called the
{\bf linearization}\index{linearization} of $f$ at the point $x$. If
the $\ssc$-continuous map $f:U\subset E\to F$ is of class $\ssc^1$,
then its tangent map
$$Tf:TU\to TF$$
is an $\ssc^0$-map. If now $Tf$ is of class $\ssc^1$, then $f:U\to
F$ is called of {\bf class} ${\mathbf {sc}}^{\boldsymbol{2}}$
\index{class!sc$^2$}.

Proceeding inductively, the map $f:U\to F$ is called of {\bf class}
$\mathbf{sc}^{\boldsymbol{k}}$ \index{class!sc$^k$} if the
$\ssc^0$-map $T^{k-1}f:T^{k-1}U\to T^{k-1}F$ is of class $\ssc^1$.
Its tangent map $T(T^{k-1}f)$ is then denoted by $T^kf$. It is an
$\ssc^0$-map $T^kU\to T^kF$.
 A map  which is of class $\ssc^k$ for every $k$ is called {\bf sc-smooth} \index{sc-smooth map} or of {\bf class}
 $\mathbf{sc}^{\boldsymbol{\infty}}$ \index{class!sc$^{\infty}$}.

Here are two useful observations which are proved in
\cite{HWZ-polyfolds1}.
\begin{prop}\label{conseq1}
If $f:U\subset E\to F$ is of class $\ssc^k$, then
$$f:U_{m+k}\to F_m$$
is of class $C^k$ for every $m\geq 0$.
\end{prop}
If we denote the usual derivative of a map $f$ by $df$ we have for
$x,h\in U_{m+1}$ the equality $df(x)(h)=Df(x)h$. In fact $Df(x)$ can
be viewed as the (unique) continuous extension of
$df(x):E_{m+1}\rightarrow F_m$ to an operator $E\rightarrow F$,
which satisfies $Df(x)(E_m)\subset F_m$  for every $m\geq 0$ and
induces continuous operators on these levels. It exists for every
$x\in U_1$ due to the definition of the class  $\ssc^1$.

 The second result is the following.

\begin{prop}\label{prop2.11}
Let $E$ and $F$ be sc-Banach spaces and let $U\subset E$ be open.
Assume that the map $f:U\rightarrow F$ is sc$^0$ and that the
induced map $f:U_{m+k}\rightarrow F_m$ is $C^{k+1}$ for every
$m,k\geq 0$. Then $f:U\rightarrow E$ is sc-smooth.
\end{prop}

Reflecting on the notion of class $\ssc^1$ one could expect for the
composition $g\circ f$ of  two such maps, that the target level
would have to drop by $2$ in order to obtain a $C^1$-map. In view of
Proposition \ref{conseq1} one might think that therefore the
composition needs not to be of class $\ssc^1$. However, this is not
the case, as the next result, the important chain rule shows.

\begin{thm}[{\bf Chain  Rule}]\label{sccomp}\index{chain rule}
Assume that $E$, $F$ and $G$ are sc-smooth Banach spaces and
$U\subset E$ and $V\subset F$  are open sets. Assume that
$f:U\rightarrow F$, $g:V\rightarrow G$ are of class $\ssc^1$ and
$f(U)\subset V$. Then the composition $g\circ f:U\to G$ is  of class
$\ssc^1$  and  the tangent maps satisfy
$$T(g\circ f)=Tg\circ Tf.$$
\end{thm}
\begin{proof}
We shall verify the properties (1) and (2) in Definition \ref{sscc1}
for $g\circ f$. The functions $g:V_1\to G$ and $f:U_1\to F$ are of
class $C^1$. Moreover, $Dg(f(x))\circ Df(x)\in {\mathcal L}(E, G)$
if $x\in U_1$. Fix $x\in U_1$ and choose
 $h\in E_1$ sufficiently small so that $f(x+h)\in V_1$.
Then, using  the postulated properties of $f$ and $g$,
\begin{equation}
\begin{split}
&g(f(x+h))-g(f(x))-Dg(f(x))\circ Df(x)h\\
&=\int_{0}^{1} Dg(tf(x+h)+(1-t)f(x))\ [f(x+h)-f(x)-Df(x)h] dt\\
&\phantom{=}+\int_{0}^{1}\bigl( [Dg(tf(x+h)+(1-t)f(x))-Dg(f(x))]
\circ Df(x)h \bigr)dt.
\end{split}
\end{equation}
Divide the first integral by the norm $\norm{h}_1$, then
\begin{equation}\label{integral1}
\begin{gathered}
\frac{1}{\norm{h}_1}\int_{0}^{1} Dg(tf(x+h)+(1-t)f(x))[f(x+h)-f(x)-Df(x)h]dt\\
= \int_{0}^{1} Dg(tf(x+h)+(1-t)f(x))\cdot \frac{1}{\norm{h}_1}
[f(x+h)-f(x)-Df(x)h]dt.
\end{gathered}
\end{equation}

If $h\in E_1$,   the maps $ [0,1]\rightarrow F_1$ defined by
$t\rightarrow tf(x+h)+(1-t)f(x) $ are  continuous and converge in
$C^0([0,1],F_1)$ to the constant map $t\rightarrow f(x)$ as
$\norm{h}_1\rightarrow 0$. Moreover,
 since $f$ is  of class $\ssc^1$,
 $$
a(h):=\frac{1}{\norm{h}_1}\bigl[ f(x+h)-f(x)-Df(x)h\bigr]$$
converges  to $0$ in $F_0$  as $\norm{h}_1\to 0$. Therefore, by the
continuity assumption (2) in Definition \ref{sscc1},
 the map $$
(t,h)\rightarrow Dg(tf(x+h)+(1-t)f(x)) [a(h)]
$$
as a map from $[0,1]\times E_1$ into $G_{0}$ converges to $0$ as
$h\rightarrow 0$,  uniformly in $t$. Therefore,   the expression in
\eqref{integral1}  converges to $0$ in $G_0$ as $h\to 0$ in $E_1$.
 Next
consider the integral
\begin{equation}\label{integral2}
\int_{0}^{1} \bigl[  Dg(tf(x+h)+(1-t)f(x))-Dg(f(x))\bigr] \circ
Df(x) \frac{h}{\norm{h}_{1}}\ dt.
\end{equation}
In view of Definition \ref{scstructure} the set of all
$\frac{h}{\norm{h}_{1}}\in E_1$ has  a compact closure in $E_{0}$.
Therefore, since $Df(x)\in {\mathcal L}(E_{0}, F_{0})$ is a
continuous map by Definition \ref{sscc1}, the closure of the set of
all
$$
Df(x) \frac{h}{\norm{h}_{1}}
$$
is compact  in $F_{0}$. Consequently, again by Definition
\ref{scstructure}, every sequence $h_n$ converging to $0$ in $E_1$
possesses a subsequence having the property that the integrand of
the integral in \eqref{integral2} converges to $0$ in $G_{0}$
uniformly in $t$. Hence the integral \eqref{integral2} also
converges to $0$ in $G_0$ as $h\to 0$ in $E_1$.  We have proved that
$$
\frac{1}{\norm{h}_{1}}\norm{g(f(x+h))-g(f(x))-Dg(f(x))\circ
Df(x)h}_{0}\rightarrow 0
$$
as $h\rightarrow 0$ in $E_1$.   Consequently,  condition (1) of
Definition \ref{sscc1} is satisfied for the composition $g\circ f$
with the linear operator
$$D(g\circ f)(x)=Dg(f(x))\circ Df(x)\in {\mathcal L}(E_0, G_0),
$$
where $x\in U_1$. We conclude that the tangent map $T(g\circ
f):TU\to TG$,
$$(x, h)\mapsto (\ g\circ f(x), D(g\circ f)(x)h\ )$$
is $\ssc$-continuous and, moreover, $T(g\circ f)=Tg\circ Tf$. The
proof of Theorem  \ref{sccomp} is complete.
\end{proof}
\mbox{}\\[2pt]
The reader should realize  that in the previous proof all conditions
on sc$^1$ maps have been used, i.e.~it just works. From Theorem
\ref{sccomp} one concludes by induction that the composition of two
$\ssc^{\infty}$-maps is also of class $\ssc^{\infty}$ and,  for
every $k\geq 1$,

\begin{equation*}
T^k(g\circ f)=T^kg\circ T^kf.\\[4pt]
\end{equation*}

An   {\bf sc-diffeomorphism} \index{sc-diffeomorphism}
  $f:U\rightarrow V$, between open subsets $U$ and $V$
   of sc-spaces $E$ and $F$ equipped with the induced sc-structure, is by
definition a homeomorphism $U\rightarrow V$ so that $f$ and $f^{-1}$
are sc-smooth.

The following remark is a continuation of Remark \ref{remake}.
\begin{rem}\label{axx}
There are other possibilities for  defining new concepts of
smoothness. For example,  we can drop the requirement of compactness
of the embedding operator $E_n\rightarrow E_m$ for $n>m$. Then it is
necessary to change the definition of smoothness in order to get the
chain rule. One needs to replace the second condition in the
definition of being sc$^1$ by the requirement that $Df(x)$ induces a
continuous linear operator $Df(x):E_{m-1}\rightarrow F_{m-1}$ for
$x\in U_{m}$ and that the map $Df:U_{m}\rightarrow {\mathcal
L}(E_{m-1},F_{m-1})$ for $m\geq 1$ is continuous. For this theory
the sc-smooth structure on $E$ given by $E_m=E$ recovers the usual
$C^k$-theory. However,  this modified theory is not  applicable to
the Gromov-Witten theory, the  Floer theory, and the SFT.
\end{rem}

\subsection{Sc-Manifolds}
Using the results so far, we can define sc-manifolds. This concept
will not yet be sufficient to describe the spaces arising in the
SFT.

\begin{defn}
Let $X$ be a second countable Hausdorff space. An {\bf sc-chart}
\index{sc-chart} of $X$ consists of a triple $(U,\varphi,E)$, where
$U$ is an open subset of $X$, $E$ a Banach space with an  sc-smooth
structure and $\varphi:U\rightarrow E$ is a homeomorphism onto an
open subset $V$ of $E$. Two such charts are sc-smoothly compatible
provided the transition maps  are  sc-smooth. An {\bf sc-smooth
atlas} \index{sc-smooth atlas} consists of a family of charts whose
domains cover $X$ so that any two charts are sc-smoothly compatible.
A maximal sc-smooth atlas is called an {\bf sc-smooth structure} on
$X$. \index{sc-smooth manifold structure} The space $X$ equipped
with a maximal sc-smooth atlas is called an {\bf
sc-manifold}.\index{sc-manifold}
\end{defn}
Let us observe that a second countable Hausdorff space which admits
an  sc-smooth atlas is metrizable and paracompact since it is
locally homeomorphic to open subsets of Banach spaces.

Assume that  the space $X$ has an  sc-smooth structure. Then it
possesses  the  filtration
 $X_m$ for all $m\geq 0$ which is induced from the filtration of the charts.
Moreover, every $X_m$ inherits the sc-smooth structure
 $(X_m)_k=X_{m+k}$ for all $k\geq 0$, denoted by $X^m$.

Next we  shall  define the {\bf tangent bundle} $p:TX\rightarrow
X^1$ in a natural way so that the tangent projection $p$ is
sc-smooth. In order to do so, we use a modification of the
definition found,  for example,  in Lang's book \cite{LA}. Namely,
consider multiplets $(U,\varphi,E,x,h)$ where $(U,\varphi,E)$ is an
sc-smooth chart, $x\in U_1$ and $h\in E$.  Call two such tuples
equivalent  if $x=x'$ and $D(\varphi'\circ\varphi^{-1})(\varphi(x))h
=h'$. An equivalence class $[U,\varphi,E,x,h]$ is called a {\bf
tangent vector} at the point $x\in X_1$. The collection of all
tangent vectors of $X$ is denoted by  $TX$. The  canonical
projection is denoted by $p:TX\rightarrow X^1$. If $U\subset X$ is
open  we introduce the subset $TU\subset TX$ by $TU=p^{-1}(U\cap
X^1)$.  For a chart $(U,\varphi,E)$ we  introduce   the map
$$
T\varphi:TU\rightarrow E^1\oplus E
$$
defined by
$$
T\varphi([U,\varphi,E,x,h]) =(x,h).
$$
One easily checks that the collection  of all triples $(TU,
T\varphi, E^1\oplus E)$ defines an sc-smooth atlas for $TX$ for
which the projection $p:TX\to X^1$ is an sc-smooth map.  The {\bf
tangent space} $T_xX$ at $x\in X_1$ is the set of equivalence
classes
$$T_xX=\{[U,\varphi, E, x,h]\vert\, h \in E\}$$
which inherits from $E$ the structure of a Banach space.
 If $x\in X_{m+1}$, then $T_xX$ possesses the
  partial filtration inherited from $E_k$, $0\leq k\leq m$.
  In particular, if $x\in X_{\infty}$, then $T_xX$ possesses
  an sc-smooth structure.

There is another class of bundles which can be defined in the
present context. These are the so-called strong (vector) bundles.
They may be viewed as a special case of strong M-polyfold bundles
which will be introduced in Section \ref{chap4}. For this reason we
shall not introduce them separately and refer the reader to Remark
\ref{specialcase}.

\section{ Splicing-Based Differential
Geometry} In this section we introduce a ``splicing-based
differential geometry''.  The fundamental concepts are splicings and
splicing cores. The splicing cores have locally varying dimensions
but admit at the same time tangent spaces. Open subsets of the
splicing cores will serve as the local models of the new global
spaces called M-polyfolds. The letter M should remind of a manifold
type space obtained by gluing together in an sc-smooth way the local
models. The M-polyfolds are equipped with substitutes for tangent
bundles, so that one is able to linearize sc-smooth maps between
M-polyfolds. In one of the follow-up papers  we go further and
introduce the notion of a polyfold which is a generalization of an
orbifold and which is on the  level of generalization needed for our
applications.

\subsection{Quadrants and Splicings}\label{sec2.1} Let us call a
subset $C$ of an sc-Banach space $W$ a {\bf partial quadrant}
\index{partial quadrant} if there is a an sc-Banach space $Q$ and a
linear sc-isomorphism $T:W\rightarrow {\mathbb R}^n\oplus Q$ mapping
$C$ onto $[0,\infty)^n\oplus Q$. If $Q=\{0\}$, then $C$ is called a
{\bf quadrant}. Observe that if $C$ and $C'$ are partial quadrants
so is $C\oplus C'$.

\begin{defn}\label{def2.1}
Assume $V$ is  an open subset of a partial quadrant $C\subset W$.
Let $E$ be an sc-Banach space and let $\pi_v:E\rightarrow E$ with
$v\in V$, be a family of projections (i.e.~$\pi_v\in {\mathcal
L}(E)$ and $\pi_v\circ \pi_v=\pi_v$) so that the induced map
\begin{gather*}
\Phi:V\oplus E\rightarrow E\\
\Phi (v, e)=\pi_v(e)
\end{gather*}
 is sc-smooth. Then the triple ${\mathcal S}=(\pi,E,V)$ is called
an {\bf sc-smooth splicing}\index{sc-smooth
splicing}\index{splicing}.
\end{defn}

The extension of the sc-smoothness definition of a map
$f:V\rightarrow F$ from an open subset of a sc-Banach space to
relatively open subsets $V\subset C\subset W$ of a partial quadrant
$C$ in an sc-Banach space, which is used in Definition \ref{def2.1},
is straight forward. One first observes that the sc-structure of $W$
induces a filtration on $V$. Now, the notion of the map $f$ to be an
$\ssc^0$-map is well-defined. Then one defines an $\ssc^0$-map
$f:V\rightarrow F$ to be of class $\ssc^1$ as in the Definition
\ref{sscc1} by replacing in there $U_1$ by $V_1$ and requiring the
existence of the limit for $x\in V_1$ and all $h\in W_1$ satisfying
$x+h\in V_1$, with a linear map $Df(x)\in {\mathcal L}(W_0,F_0)$.
Moreover, the tangent bundle $TV$ of the set $V$ is defined as usual
by $TV=V^1\oplus W$ together with the $\ssc^0$-projection map
$TV\rightarrow V^1$.

 Every splicing ${\mathcal S}=(\pi, E, V)$ is accompanied by the
{\bf complementary splicing} ${\mathcal S}^c=(1-\pi, E, V)$ where
$1-\pi$ stands for the family of projections
$(1-\pi_v)(e)=e-\pi_v(e)$ for $(v, e)\in V\oplus E$. This way the
splicing decomposes the set $V\oplus E$ naturally into a fibered sum
over the parameter set $V$. Indeed, $(v, e)\in V\oplus E$ can be
decomposed as
$$
(v, e)=(v, e_v+e_v^c)
$$
where  $\pi_v (e)=e_v$ and $(1-\pi_v)(e)=e_v^c$. The  splicing cores
$K^{\mathcal S}$ and $K^{{\mathcal S}^c}$ can be viewed as bundles
over $V$ (with linear Banach space fibers, which however change
dimensions). Their Whitney sum over $V$
$$
K^{\mathcal S}\oplus_V K^{{\mathcal S}^c}=\{(v,a,b)\in V\oplus
E\oplus E\ | \pi_v(a)=a,\ \pi_v(b)=0\}
$$
is naturally diffeomorphic to $V\oplus E$. The name splicing comes
from the fact that it defines a decomposition of $V\oplus
E\rightarrow V$, by `splicing' it along $V$.

We should point out that the sc-smoothness of the mapping $(v,
e)\mapsto \pi_v (e)$ is a rather weak requirement allowing the
dimension of the images of the projections $\pi_v $ to vary locally
in the parameter $v\in V$. The reader can find illustrations and
examples in \cite{HWZ-polyfolds1}.

Since $\pi_v$ is a projection,
\begin{equation}\label{nnn}
\Phi(v,\Phi(v,e))=\Phi(v,e).
\end{equation}
The left-hand side is the composition of $\Phi$ with the sc-smooth
map $(v,e)\rightarrow (v,\Phi(v,e))$.  For  fixed $(v,\delta v)\in
TV$ we introduce the map
\begin{equation}\label{equation2.2}
\begin{gathered}
P_{(v,\delta v)}:TE\rightarrow TE\\
(e,\delta e)\rightarrow (\Phi(v,e),D\Phi (v,e)(\delta v,\delta e)).
\end{gathered}
\end{equation}
It has  the property that the induced map
$$
TV\oplus TE\rightarrow TE:(a,b)\rightarrow P_a(b)
$$
is sc-smooth because,   modulo the identification $TV\oplus
TE=T(V\oplus E)$,   it is equal to the tangent map of $\Phi$.  From
\eqref{nnn} one obtains by means of the chain rule (Theorem
\ref{sccomp}) at the points $(v, e)\in (V\oplus E)_1$, the formula

$$
D\Phi(v,\Phi(v,e))(\delta v, D\Phi(v,e)(\delta v,\delta
e))=D\Phi(v,e)(\delta v,\delta e)
$$
\mbox{}\\[1pt]
and,  together with the definition of $P$,  one computes\mbox{}\\[1pt]
\begin{equation*}
\begin{split}
 P_{(v,\delta v)}\circ P_{(v,\delta v)}(e,\delta e)&
=P_{(v,\delta v)}(\pi_v(e),D\Phi(v,e)(\delta v,\delta e))\\
&=(\pi_v^2(e),D\Phi (v,\pi_v(e))(\delta v, D\Phi(v,e)(\delta
v,\delta
e)) )\\
&=(\pi_v(e),D\Phi(v,e)(\delta v,\delta e)) = P_{(v,\delta
v)}(e,\delta e).
\end{split}
\end{equation*}
\mbox{}\\[1pt]
Consequently, $P_{(v,\delta v)}$ is a projection, which of course
can be identified with the tangent $T\pi$ of the map $\pi:V\oplus E\to E$, defined by
$\pi (v, e)=\pi_v (e)$, via
$$
P_{(v,\delta v)}(e,\delta e)= T\pi ((v,e),(\delta v,\delta e)).
$$
In the following we shall write $(T\pi)_{(v,\delta v)}$ instead of
$P_{(v,\delta v)}$.  Hence the triple
$$T{\mathcal S}=(T\pi,TE,TV)$$
is an  sc-smooth splicing,  called  the {\bf  tangent splicing of}
${\mathcal S}$.\index{tangent splicing}
\begin{defn}[{\bf Splicing Core}]\label{splicingcore}
If  ${\mathcal S}=(\pi,E,V)$ is  an sc-smooth splicing, then the
associated  {\bf splicing core} \index{splicing core}  is  the image
bundle of the projection $\pi$ over $V$, i.e., it is the subset
$K^{\mathcal S}\subset V\oplus E$ defined by\\[1pt]
\begin{equation}\label{scpor1}
K^{\mathcal S}:=\{ (v, e)\in V\oplus E\vert\: \pi_v(e)=e\}.
\end{equation}
\end{defn}

If the dimension of $E$ is finite, the images of the projections
$\pi_v$ have all the same rank so that the splicing core is a smooth
vector bundle over $ V$. If, however, the dimension of $E$ is
infinite, then the ranks of the fibers can change with the parameter
$v$
thanks to the definition of $\ssc$-smoothness.
This truly infinite dimensional phenomenon is crucial for our purposes.\\[4pt]

The {\bf splicing core of the tangent splicing}  $T{\mathcal S}$ is
the set
\begin{equation}\label{spcor2}
K^{T{\mathcal S}}=\{(v, \delta v, e,\delta e)\in TV\oplus TE\vert
\:\: (T\pi)_{(v,\delta v)}(e,\delta e)=(e, \delta e)\}.
\end{equation}
The mapping
\begin{equation*}
K^{T{\mathcal S}}\rightarrow {(K^{\mathcal S})}^1: (v,\delta
v,e,\delta e)\mapsto  (v,e) \in V_1\oplus E_1
\end{equation*}
is the canonical projection. The fiber over every  point  $(v,e)\in
{(K^{\mathcal S})}^1$ is a subspace $K^{T{\mathcal S}}_{(v,e)}$ of
the Banach space $W\oplus E$.  If $(v, e)$ is on level $m+1$, then
$K^{T{\mathcal S}}_{(v, e)}$ has well defined levels $k\leq m$. The
tangent splicing $K^{T{\mathcal S}}$ has well defined bi-levels $(m,
k)$ with $k\leq m$. Indeed, assuming for simplicity that
$W=\R^n\oplus Q$, then $V\subset C\subset W$ and $TV=V^1\oplus W$
and we can define for $0\leq k\leq m$,
$$
(K^{T{\mathcal S}})_{m,k}=\{(v, \delta v, e,\delta e)\in V_{m+1}
\oplus W_k\oplus E_{m+1}\oplus E_k\vert \: (T\pi)_{(v,\delta
v)}(e,\delta e)=(e,\delta e)\}.
$$
The projection $K^{T{\mathcal
S}}\rightarrow {(K^{\mathcal S})}^1: (v,\delta v,e,\delta e)\mapsto
(v,e)$
 maps level $(m,k)$ points to level $m$ points. We may view $k$
 as the fiber regularity and $m$ as the base regularity.
Note that a  point $e$ of $E^j$ of regularity $m$ has regularity
$m+j$ as a point in $E$. The following is one of our main
definitions.
\begin{defn}\label{vnewdefn2.4}
A {\bf local  M-polyfold model} \index{local m-polyfold model}
consists of a pair $(O,{\mathcal S})$ where $O$ is an open subset of
the splicing core $K^{\mathcal S}\subset V\oplus E$ associated with
the sc-smooth splicing ${\mathcal S}=(\pi,E,V)$. The {\bf tangent of
the local M-polyfold model} \index{tangent of the local M-polyfold
model}  $(O,{\mathcal S})$ is the object defined by
$$
T(O,{\mathcal S}) = (K^{T{\mathcal S}}|O^1,T{\mathcal S})
$$
where $K^{T{\mathcal S}}|O^1$ denotes the collection of all points
in $K^{T{\mathcal S}}$ which project under the canonical projection
$ K^{T{\mathcal S}}\rightarrow {(K^{\mathcal S})}^1$ onto  the points in $O^1$.
\end{defn}
There is  the  natural projection
$$
K^{T{\mathcal S}}\vert O^1\to O^1:(v,\delta v,e,\delta e)\rightarrow
(v,e).
$$
In the following we shall simply write $O$ instead of $(O,{\mathcal
S})$, but keep in mind that ${\mathcal S}$ is part of the structure.
With this notation the tangent $TO=T(O, {\mathcal S})$ of the open
subset $O$ of the splicing core $K^{\mathcal S}$ is the set

\begin{equation}\label{kts}
TO=K^{T{\mathcal S}}\vert O^1.\mbox{}\\[2pt]
\end{equation}
\noindent Note that on an open subset $O$ of a splicing core there
is an induced filtration. Hence we may talk about sc$^0$-maps. We
will see in the next section  that there is  also a well-defined
notion of a $\ssc^1$-map in this setting. We shall see in the
applications presented in
\cite{HWZ-polyfolds1,HWZ-polyfolds2,HWZ4,HWZ5} that analytical
limiting phenomena, like bubbling-off, occurring in symplectic field
theory, Gromov-Witten theory and Floer Theory are smooth within the
splicing world.

\subsection{Smooth Maps between Splicing Cores}\label{smspco}

The aim of this section is to introduce the concept of an
$\ssc^1$-map between local M-polyfold models. We will construct the
tangent functor and show the validity of the chain rule. At that
point we will have established  all the ingredients for building the
{\bf ``splicing-differential geometry"} mentioned in the
introduction.

Consider two open subsets $O\subset K^{\mathcal S}\subset V\oplus E$
and $O'\subset K^{\mathcal S'}\subset V'\oplus E'$ of splicing cores
belonging to the splicings
 ${\mathcal S}=(\pi, E,V)$ and ${\mathcal S'}=(\pi', E', V')$.
 The open subsets $V$ and $V'$ of partial quadrants are contained in the sc-Banach spaces $W$ resp. $W'$.
Consider an $\ssc^0$-map
\begin{equation*}\label{map1}
f:O\to O'.
\end{equation*}

If $O$ is an open subset of the splicing core  $K^{\mathcal
S}\subset V\oplus E$ we define the subset $\wh{O}$ of $V\oplus E$ by
$$\wh{O}=\{(v, e)\in V\oplus E\vert \: (v, \pi_v(e))\in O\}.$$
Clearly, $\wh{O}$ is open in $V\oplus E$ and can be viewed as a
bundle $\wh{O}\rightarrow O$ over $O$. This bundle will be important
in subsection \ref{section4.2}, where the crucial notion of a filler
is being introduced.

\begin{defn}\label{def2.7}
The $\ssc^0$-continuous map $f:O\to O'$ between open subsets  of
splicing cores is called of {\bf class} ${\mathbf
{sc}}^{\boldsymbol{1}}$ if the map
\begin{gather*}
\wh{f}:\wh{O}\subset V\oplus E\to W'\oplus E'\\
\wh{f}(v, e)=f(v, \pi_v(e))
\end{gather*}
is of class $\ssc^1$.
\end{defn}
According to the splitting of the image space we set
$$\wh{f}(v,e)=(\wh{f}_1(v, e), \wh{f}_2(v, e))
\in K^{{\mathcal S}'}\subset W'\oplus E'.
$$
 The {\bf tangent map} \index{tangent map between splicings} $T\wh{f}$ associated
 with the $\ssc^1$-map $\wh{f}$  is defined as
\begin{equation}\label{map2}
T\wh{f}(v,\delta v, e,\delta e):=(T\wh{f}_1(v,\delta v, e,\delta
e),T\wh{f}_2(v,\delta v, e,\delta e)).
\end{equation}
The map $T\wh{f}$ is of class $\ssc^0$.
\begin{lem}\label{tanmaplem}
The tangent map $T\wh{f}$ satisfies $T\wh{f}(K^{T{\mathcal S}}\vert
O^1)\subset K^{T{\mathcal S'}}\vert {O'}^1$ and hence induces a map
$$
K^{T{\mathcal S}}\vert O^1\to K^{T{\mathcal S'}}\vert {O'}^1
$$
which we denote by $Tf$. In the simplified notation of \eqref{kts},
we have
$$Tf:TO\to TO'.$$
\end{lem}
\begin{proof}
Denote by $\pi_{v'}'$ the family of projections associated with the
splicing ${\mathcal S'}$. Since $f:O\to O'$,  we have by definition
of  the splicing core $K^{\mathcal S'}$ the formula
$$\pi'_{\wh{f}_1(v, e)}(\wh{f}_2(v, e))=\wh{f}_2(v, e).$$
Differentiating this  identity in the variable $(v, e)$ we obtain
\begin{equation}\label{map22}
\begin{split}
D\wh{f}_2(v, e)(\delta v, \delta e)&=
D_{v'}\pi'_{\wh{f}_1(v, e)}(\wh{f}_2(v, e))\circ
D\wh{f}_1(v, e)(\delta v, \delta e)\\
&\phantom{=}+\pi'_{\wh{f}_1(v, e)}\circ D\wh{f}_2(v, e)(\delta v,
\delta e).
\end{split}
\end{equation}
Set
$$
v'=\wh{f}_1(v, e),\ e'=\wh{f}_2(v, e),\ \delta v'=D\wh{f}_1(v,
e)(\delta v, \delta e),\ \delta e'=D\wh{f}_2(v, e)(\delta v, \delta
e).
$$ Then  \eqref{map22} implies using the definition
\eqref{equation2.2} of the projection $(T\pi')_{(v', \delta v')}$
associated with the splicing $T{\mathcal S}'$ that
$$(T\pi')_{(v', \delta v')}(e', \delta e')=(e', \delta e').$$
So,  indeed $Tf(v, \delta v, e,\delta e)=(v', \delta v', e', \delta
e')\in K^{T{\mathcal S}'}$ as was claimed.
\end{proof}
Note that the order of the terms in the tangent map $Tf$ resp.
$T\hat{f}$ of an $\ssc^1$-map $f:O\rightarrow O'$ is different from
the order of terms in the classical notation. Writing $f=(f_1,f_2)$
according to the splitting of the image space into the distinguished
splicing parameter part and the standard part, the classical
notation for the tangent map would be $Tf=((f_1,f_2),(Df_1,Df_2))$
whereas our convention is $Tf=((f_1,Df_1),(f_2,Df_2))$. This rather
unorthodox ordering of the data has been chosen so that the tangent
of a splicing is again a splicing.

 The reader could work out as an
example the situation where the splicings have the constant
projection $Id$.
\begin{thm} [{\bf Chain Rule}]\label{chainrule1}
Let $O,O',O''$ be open subsets of splicing cores and let the the
maps $f:O\rightarrow O'$ and $g:O'\rightarrow O''$ be of class
$\ssc^1$. Then  the composition $g\circ f$ is also of class $\ssc^1$
and   the tangent map satisfies \\[1pt]
$$T(g\circ f) =Tg\circ Tf.$$
\mbox{}
\end{thm}
\begin{proof}
This is a consequence of the sc-chain rule (Theorem \ref{sccomp}),
the definition of the tangent map and the fact that our reordering
of the terms in our definition \eqref{map2} of the tangent map is
consistent. Indeed, from  Definition \eqref{map2} we deduce
 \begin{equation*}
 \begin{split}
&T(g\circ f)(v,\delta v,e,\delta e) = (T(\what{g}_1\circ
\what{f})(v,e,\delta v,\delta e),
T(\what{g}_2\circ\what{f})(v,e,\delta v,\delta e))\\
&\phantom{===}= ((T\hat{g}_1)\circ (T\hat{f})(v,e,\delta v,\delta
e),
(T\what{g}_2)\circ  (T\what{f})(v,e,\delta v,\delta e))\\
&\phantom{===}= (Tg)(T\what{f}_1(v,e,\delta v,\delta
e),T\what{f}_2(v,e,\delta
v,\delta e))\\
&\phantom{===}= (Tg)\circ (Tf)(v,\delta v,e,\delta e)
 \end{split}
\end{equation*}
and the proof is complete.
\end{proof}

Given an  sc$^1$-map $f:O\rightarrow O'$ between open sets of
splicing cores we obtain, in view of Lemma \ref{tanmaplem},  an
induced tangent map $Tf:TO\rightarrow TO'$. Since $TO$ and $TO'$ are
again open sets in the splicing cores $K^{T{\mathcal S}}$ and
$K^{T{\mathcal S'}}$ we can  iteratively define  the notion of $f$
to be of {\bf class $\ssc^k$} and of $f$ to be sc-smooth.

\begin{defn}\label{tangentspace1}
Let $O$ be an open subset of a splicing core $K^{\mathcal S}$ and
$(v, e)\in O_1$. The {\bf tangent space} to $O$ {\bf at the point}
$(v, e)$ is the Banach space

 \begin{equation}\label{map4}
 T_{(v, e)}O=\{(\delta v, \delta e)\in W\oplus E\, \vert\: \: (v, \delta v, e,\delta e)\in TO\}.
 \end{equation}
 \end{defn}

We then have  $$TO=\bigcup_{(v, e)\in O_1}T_{(v, e)}O.$$
  If $f:O\to O'$ is a homeomorphism so that $f$ and $f^{-1}$ are sc-smooth,
   our tangent map $Tf$ defined in \eqref{map2} induces the linear isomorphism
 $$
Tf(v,e):T_{(v,e)}O\rightarrow T_{f(v,e)}O'.
$$
We recall from Section \ref{sec2.1} that the space
 $TO$ has a bi-filtration $(TO)_{(m,k)}$ for $0\leq k\leq m$,
 so that the natural projection
 $$
 TO\rightarrow O^1
 $$
 maps level $(m,k)$ points to level $m$ points and is $\ssc$-smooth.
 The projection map $TO\rightarrow O^1$ is sc-smooth.

\subsection{M-Polyfolds}\label{section2.4}
 Now we are  able to
introduce the notion of  an M-polyfold. The ``M" indicates the
``manifold flavor" of the polyfold. A general polyfold will be a
generalization of an orbifold.
\begin{defn}\label{DEF38}
Let $X$ be a second countable Hausdorff space. An {\bf M-polyfold
chart} \index{m-polyfold chart}for $X$ is a triple
$(U,\varphi,{\mathcal S})$, in which $U$ is an open subset of $X$,
${\mathcal S}=(\pi,E,V)$ an sc-smooth splicing and
$\varphi:U\rightarrow K^{\mathcal S}$ a homeomorphism onto an open
subset $O$ of the splicing core $K^{\mathcal S}$ of ${\mathcal S}$.  Two
charts are called compatible if the transition maps between open
subsets of splicing cores are sc-smooth in the sense of Definition
\ref{def2.7}. A maximal atlas of sc-smoothly compatible M-polyfold
charts is called an  {\bf M-polyfold structure}\index{m-polyfold
structure} on $X$.
\end{defn}

An  M-polyfold is necessarily metrizable by  an argument similar to
the one  used already for sc-manifolds. Each splicing core
$K^{\mathcal S}$ carries the structure of an M-polyfold with the
global chart being the identity.

The concept of a map $f:X\to X'$ between M-polyfolds to be of class
$\ssc^{0}$ or $\ssc^k$ or to be sc-smooth is, as usual, defined by
means of local charts.
\begin{defn}\label{localcharts}
The mapping $f:X\to X'$ between two M-polyfolds is called of {\bf
class $\ssc^{\mathbf 0}$} resp. {\bf $\ssc^{\mathbf k}$} or called
{\bf sc-smooth} if for every point $x\in X$ there exists  a chart
$(U, \varphi, {\mathcal S})$ around $x$  and a chart $(U', \varphi',
{\mathcal S}')$ around $f(x)$ so that $f(U)\subset U'$ and
$$\varphi'\circ f\circ \varphi^{-1}:\varphi (U)\to \varphi'(U')$$
is of class $\ssc^0$ resp. $\ssc^k$ or sc-smooth.
\end{defn}

 In order to define the {\bf tangent space}
 $T_xX$ of the $M$-polyfold $X$ at the point $x\in X_1$,
 we proceed as in the case of sc-manifolds in section \ref{section2.4}.
This time we consider equivalence classes  of multiplets $(U, \varphi, S, x, h)$ in which $(U, \varphi, S)$ is an $M$-polyfold chart, $x$ is a point in $U_1$ and $h\in T_{\varphi (x)}O$, where $O=\varphi (U)\subset K^{{\mathcal S}}$ is the open set of the splicing core. The above multiplet is equivalent to $(U', \varphi', S', x', h')$ if $x=x'$ and if $T(\varphi'\circ \varphi^{-1})(\varphi (x))h=h'$, where the tangent map
$$
T(\varphi'\circ \varphi^{-1})(\varphi (x)):T_{\varphi (x)}O\to
T_{\varphi'(x)}O'
$$
is defined in section \ref{smspco}. The tangent space
is now defined as the set of equivalence classes
$$T_xX=\{[U, \varphi, S, x, h]\vert \, h\in T_{\varphi (x)}O\}.$$
It inherits the structure of a Banach space from the tangent
 space $T_{\varphi (x)}O$. If $x\in X_{m+1}$, then $T_xX$ possesses
  a partial filtration for $0\leq k\leq m$ induced from the
   partial filtration of $T_{\varphi (x)}O$.
   The tangent space at a smooth point $x\in X_{\infty}$
   possesses an sc-smooth structure.

Let now $f:X\to X'$ be a map between $M$-polyfolds of class $\ssc^k$ for $k\geq 1$. In two $M$-polyfold charts $(U, \varphi, S)$ and $(U', \varphi', S')$ around the points $x\in U_1$ and $f(x)\in U_1'$, the map $f$ is represented  by the $\ssc^k$-map $\psi:\varphi'\circ f\circ \varphi^{-1}:O\to O'$ between open sets of splicing cores. The tangent map $T\psi (x):T_{\varphi (x)}O\to T_{\psi (x)}O'$ defines a unique continuous linear map
$$T_xf:T_xX\to T_{f(x)}X'$$
between the tangent spaces, called the {\bf tangent map} of $f$ at the point $x$, mapping the equivalence class $[U, \varphi, S, x,h]$ into the class $[U', \varphi', S', f(x), h']$ in which
$$h'=T(\varphi'\circ f\circ \varphi^{-1})(\varphi (x))\cdot h.$$
If $x$ is a smooth point of $X$ and if $f$ is an sc-smooth map, then the tangent map $T_xf$ is an sc-operator as defined in section \ref{section2.2}.

Let us note the following useful  result about sc-smooth partitions
of unity.
\begin{thm}\label{partitionunity}
Let $X$ be an  M-polyfold with local models being splicing cores
build on sc-Hilbert spaces (An sc-Hilbert space consists of a
Hilbert space equipped with an sc-structure. It is not required that
the Banach spaces $E_m$ for $m\geq 1$ are Hilbert spaces.). Assume
that $(U_{\lambda})_{\lambda\in\Lambda}$ is an open covering of $X$.
Then there exists a subordinate sc-smooth partition of unity
$(\beta_{\lambda})_{ \lambda \in \Lambda}$.
\end{thm}
The statement follows along the lines of a proof for Hilbert
manifolds in \cite{LA}. The product $X\times Y$ of two M-polyfolds
is in a natural way an M-polyfold. Indeed, if $(U,\varphi,{\mathcal
S})$ and $(W,\psi,{\mathcal T})$ are $M$-polyfold charts for  $X$
and $Y$ respectively,  one obtains the product chart  $(U\times
W,\varphi\times\psi,{\mathcal S}\times {\mathcal T})$ for $X\times
Y$, with the {\bf product splicing}\index{product splicing}
\begin{equation*}
\begin{split}
{\mathcal S}\times{\mathcal
T}&=(\pi,E,V)\times(\rho,F,V')\\
&=(\sigma,E\oplus F,V\oplus V')
\end{split}
\end{equation*}
where  $\sigma_{(v,v')}=\pi_v\oplus\rho_{v'}$ is the family of
projections. There are several possible notions of sub-polyfolds (we
suppress the $M$ in the notation). We shall describe one of them in
Section \ref{immsub} below and refer the reader to
\cite{HWZ-polyfolds1} for a comprehensive treatment.

\subsection{Corners and Boundary Points}\label{section2.5}

In this section we will  prove  the extremely important fact that
sc-smooth maps  are able to recognize  corners. This will be crucial
for the SFT because  most of its  algebraic structure  is a
consequence of the corner structure.

Let $X$ be a M-polyfold. Around a point $x\in X$ we take a
M-polyfold chart $\varphi:U\rightarrow K^{\mathcal S}$ where
$K^{\mathcal S}$ is the splicing core associated with the splicing
${\mathcal S}=(\pi,E,V)$. Here $V$ is an open subset of a partial
quadrant $C$ contained in the sc-Banach space $W$. By definition
there exists a linear isomorphism from $W$ to ${\mathbb R}^n\oplus
Q$ mapping  $C$ onto $[0,\infty)^n\oplus Q$.  Identifying the
partial quadrant $C$ with $[0,\infty )^n\oplus Q$ we shall use the
notation $\varphi=(\varphi_1,\varphi_2)\in [0,\infty)^n \oplus
(Q\oplus E)$ according to the splitting of the target space of
$\varphi$. We associate with the point $x\in U$ the integer  $d(x)$
defined by
\begin{equation}\label{di}
d(x)=\sharp \{\text{coordinates of $\varphi_1(x)$ which are equal to
$0$}\}.
\end{equation}

\begin{thm}\label{cornerthm}
The map $d:X\rightarrow {\mathbb N}$ is well-defined and does not
depend on the choice of  the M-polyfold chart $\varphi:U\to
K^{{\mathcal S}}$. Moreover, every point $x\in X$ has  an open
neighborhood $U'$ satisfying
$$
d(y)\leq d(x)\ \hbox{for all}\ \ y\in U'.
$$
\end{thm}

\begin{defn}\label{degindex}
The map $d:X\to \N$ is called the {\bf degeneracy index} of
$X$.\index{degeneracy index of the M-polyfold}
\end{defn}

The map $d$ will play an important role in our Fredholm theory with
operations presented in \cite{HWZ3.5}. A point  $x\in X$ satisfying
$d(x)=0$ is called an {\bf interior point}. A point satisfying
$d(x)=1$ is called a {\bf good boundary point}\index{good boundary
point}. A point with $d(x)\geq 2$ is called a {\bf  corner}. In
general,  the integer $d(x)$ is the order \index{order of the
corner} of the corner. \index{corner}

\begin{proof}[Proof of Theorem \ref{cornerthm}]
Consider two $M$-polyfold charts $\varphi:\wh{U}\subset X\to K^{S}$
and $\varphi':\wh{U}'\subset X \to K^{{\mathcal S}'}$ such that
$x\in \wh{U}\cap \wh{U}'$. Introducing   the open subsets $U=\varphi
(\wh{U}\cap \wh{U}')$ and $U'=\varphi'(\wh{U}\cap \wh{U}')$ of
$K^{{\mathcal S}}$and $K^{\mathcal S'}$ resp., and setting $\varphi
(x)=(r, a)$ and $\varphi' (x)=(r', a')$ we define the
sc-diffeomorphism $\Phi :U\to U'$  by $\Phi=\varphi'\circ
\varphi^{-1}$. Obviously, $\Phi (r, a)=(r', a')$. Now the proof of
Theorem \ref{cornerthm} reduces to the following proposition.

\begin{prop} \label{cornerprop} Let ${\mathcal S}=(\pi,E,V)$ and
${\mathcal S}'=(\pi', E',  V')$ be two   splicings having the
parameter sets $V=[0,\infty)^k\oplus Q$ and
$V'=[0,\infty)^{k'}\oplus Q'$. Assume that   $U$ and $U'$ are open
subsets of the splicing cores $K^{{\mathcal S}}$ and $K^{{\mathcal
S'}}$ containing the points $(r, a)$ and $(r', a')$ with $r\in
[0,\infty)^k$ and $r'\in [0,\infty)^{k'}$ and assume that the map
$$
\Phi:U\rightarrow U'
$$
is an sc-diffeomorphism mapping $(r,a)$ to $(r',a')$. Then $r$ and
$r'$ have the same number of vanishing coordinates.
\end{prop}
\begin{proof}
We first prove the assertion  under the additional assumption that
the point  $p_0=(r, a)$ belongs to $ U_{\infty}$. Then  the image
point $q_0=(r',a')=\Phi (p_0)$ belongs to $U'_{\infty}$. Denote  by
$J$ the subset of $\{1,\cdots ,k\}$ consisting of those indices $j$
for which $r_j=0$. Similarly, $j'\in J'\subset \{1,\cdots ,k'\}$ if
$r'_j=0$. Denoting  by $\sharp r$ and $\sharp r'$ the cardinalities
of $J$ and $J'$  we claim that $\sharp r=\sharp r'$.  Since $\Phi$
is an sc-diffeomorphism it suffices to prove the inequality $\sharp
r\geq \sharp r'$ since the inequality has to also  hold  true for
the sc-diffeomorphism $\Phi^{-1}$. Write $a=(q,e)$. If
$\pi_{(r,q)}(e)=e$, then differentiating $\pi_{(r,q)}\circ
\pi_{(r,q)} (e)=\pi_{(r,q)} (e)$ in $(r,q)$ one finds $\pi_{(r,q)}
\circ D_{(r,q)}(\pi_{(r,q)}(e))=0$ so that
$D_{(r,q)}(\pi_r(e))(\delta r,\delta q)$ is contained in the range
of $\text{id}-\pi_{(r,q)}.$ Therefore, given $(r,a)\in U_{\infty}$
satisfying $\pi_{(r,q)}(e)=e$ and given $\delta r\in \R^k$, and
$\delta q\in Q_\infty$, there exists $\delta e\in E_{\infty}$
solving
\begin{equation}\label{corner1}
\delta e=\pi_{(r,q)} (\delta e)+D_{(r,q)} (\pi_{(r,q)} (e))[(\delta
r,\delta q)].
\end{equation}
In particular, taking $\delta r\in \R^k$ with $(\delta r)_j=0$ for
$j\in J$, and a smooth $\delta q$, there exists $\delta e\in
E_{\infty}$ solving the equation \eqref{corner1}. This is equivalent
to $((\delta r,\delta q), \delta e)\in  (T_{(r,a)}U)_{\infty}$.
Introduce the path
$$
\tau \mapsto p_{\tau}=(r+\tau \delta r,q+\tau\delta q, \pi_{(r+\tau
\delta r,q+\tau\delta q)}(e+\tau \delta e))
$$
for $\abs{\tau }<\rho$ and  $\rho$ small. From  $(r, a)\in
U_{\infty}$ and $\delta e\in E_{\infty}$ one concludes $p_{\tau}\in
U_{\infty}$. Moreover,  considering $\tau \to p_{\tau}$ as a map
into $U_m$ for  $ m\geq 0$, its derivative at $\tau=0$ is equal to
$(\delta r,\delta q, \delta e)$. Fix a level $m\geq 1$  and consider
for $\rho>0$ sufficiently small the map
$$
(-\rho,\rho)\rightarrow {\mathbb R}^{k'}\oplus Q'_m\oplus E_{m}':
\tau\rightarrow \Phi(p_{\tau}).
$$
The map $\Phi:U\to U'$ is $C^1$ as a map from $U_{m+1}\subset
\R^k\oplus Q_{m+1}\oplus  E_{m+1}$ into $\R^{k'}\oplus Q_m'\oplus
E_m'$. Its derivative $d\Phi (r,q, e):\R^k\oplus Q_{m+1}\oplus
E_{m+1}\to \R^{k'}\oplus Q'_m \oplus E_m'$ has an extension to the
continuous linear operator $D\Phi (r,q, e):\R^k\oplus Q_m\oplus
E_m\to \R^{k'}\oplus Q'_m\oplus E_m'$. Since $\Phi$ is a
sc-diffeomorphism the extension $D\Phi(r,q, e):\R^k\oplus Q_m\oplus
E_m\to \R^{k'}\oplus Q_m'\oplus E_m'$ is a bijection. Thus, since
$\delta q\in Q_\infty$ and $\delta e\in E_{\infty}$,
\begin{equation}\label{corner2}
\begin{split}
\Phi (p_{\tau})&=\Phi (p_0)+\tau \cdot d\Phi (p_0)[\delta r,
\delta q,\delta e]+o_{m}(\tau)\\
&=q_0+\tau \cdot D\Phi (p_0)[\delta r, \delta q, \delta
e]+o_{m}(\tau)
\end{split}
\end{equation}
where $o_{m}(\tau)$ is a function taking values in $\R^{k'}\oplus
Q_k'\oplus  E_m'$ and satisfying $\frac{1}{\tau} o_{m}(\tau )\to
0\quad \text{as $\tau \to 0$}.$ Introduce the sc-continuous linear
functionals $\lambda_{j'}:\R^{k'}\oplus Q'\oplus  E'\to \R$ by
$$\lambda_{j'}(s',q',h')=s_{j'}'
$$
where $j'\in \{ 1,\ldots ,k'\}$. Then
$$
\lambda_{j'}\circ \Phi(p_{\tau})\geq 0
$$
for $\abs{\tau}<\rho$ and for $j'\in \{1,...,k'\}$. Applying for
$j'\in J'$ the functional $\lambda_{j'}$ to both sides of
\eqref{corner2} and using that for $j'\in J'$ we have
$\lambda_{j'}(\Phi(p_0))=\lambda_{j'}(q_0)=0$ we conclude for
$\tau>0$
\begin{equation*}\label{corner3}
\begin{split}
0\leq \frac{1}{\tau}\cdot \lambda_{j'} \bigl[ \Phi (p_{\tau})\bigr]
&=\frac{1}{\tau}\cdot \lambda_{j'}\bigl[\Phi (p_0)+
\tau \cdot D\Phi (p_0)[\delta r,\delta q, \delta e]+o_{m}(\tau) \bigr]\\
&=\lambda_{j'} \bigl[  D\Phi (p_0)[\delta r, \delta q, \delta
e]\bigr] +\lambda_{j'}\bigl(\frac{o_{m}(\tau )}{\tau}\bigr).
\end{split}
\end{equation*}
Passing to the limit $\tau\rightarrow 0^+$ we find
\begin{equation*}\label{corner4}
0\leq \lambda_{j'} (D\Phi (p_0)[\delta r, \delta q,\delta e])
\end{equation*}
and replacing $(\delta r, \delta q,\delta e)$ by $(-\delta r,-\delta
a, -\delta e)$ we obtain the equality  sign. Consequently,
\begin{equation}\label{corner10}
\lambda_{j'} (D\Phi (p_0)[\delta r, \delta q,\delta e])=0, \quad
j'\in J'
\end{equation}
for all $[\delta r, \delta q,\delta e]\in \R^k\oplus Q_\infty\oplus
E_{\infty}$ satisfying
$$
\pi_{(r,q)} (\delta e)+D_{(r,q)} (\pi_{(r,q)}(e))[(\delta r,\delta
q)]=\delta e
$$
and $(\delta r)_j=0$ for all $j\in J$. Introduce the codimension
$\sharp r$ subspace $L$ of the tangent space $T_{(r,q,e)}U\subset
K^{T{\mathcal S}}$ which we may view as a subset of $\R^k\oplus
Q_\infty\oplus E_{\infty}$ by
\begin{eqnarray*}
L&=& \{(\delta r, \delta q,\delta e)\in \R^k\oplus Q_\infty\oplus
E_{\infty}\vert\ \text{$\pi_{(r,q)}(\delta e)+D_{(r,q)} (\pi_{(r,q)}
(e))(\delta
r,\delta q)=\delta e$}\\
&& \text{and $\: (\delta r)_j=0\ $ for all $ j\in J$}\: \}.
\end{eqnarray*}
Then,   in view of  \eqref{corner10},
\begin{eqnarray*}
D\Phi(r, q,e)L  &\subset& \{(\delta  r', \delta q',\delta
e')\vert\pi_{(r',q')}'(\delta
e')+D_{(r',q')}(\pi_{(r',q')}'(e'))(\delta
r',\delta q')=\delta e'\\
&&\text{and}\ (\delta r')_{j'}=0\ \text{for all } \ j'\in J'\}.
\end{eqnarray*}
 Because the subspace on the right
hand side has codimension $\sharp r'$ in $T_{(r',q',e')}U'$ and
since $D\Phi (r,q, e)$,  being a bijection, maps $L$ onto a
codimension $\sharp r$ subspace of $T_{(r',q',e')}U'$, it follows
that $\sharp r'\leq \sharp r$, as claimed.

Next we shall prove the general case. For this we  take
$p_0=(r,q,e)$ in $U_{0}$, so that the image point
$(r',q',e')=\Phi(r,q,e)$ belongs to $U'_0$. Arguing by contradiction
we may assume without loss of generality that $\sharp r > \sharp
r'$, otherwise we replace $\Phi$ by $\Phi^{-1}$.  Since $U_{\infty}$
is dense in $U_0$ we find a sequence  $(r,q_n,e_n)\in U_{\infty}$
satisfying $\pi_{(r,q_n)} (e_n)=e_n$ and $(r,q_n,e_n)\rightarrow
(r,q,e)$ in $U_{0}$.  By the previous discussion $\sharp r=\sharp
r'_n$ where $(r_n',q_n',e_n')=\Phi(r,q_n,e_n)$.    Since $\Phi$ is
sc-smooth, we have $(r_n',q_n',e_n')\rightarrow (r',q',e')$ in
$U'_{0}$ and $\pi_{(r',q')}'(e')=e'$.    From this convergence we
deduce $\sharp r'\geq  \sharp r_n'$ so that  $\sharp r' \geq \sharp
r$ contradicting  our assumption. The proof of Proposition
\ref{cornerprop} is complete.
\end{proof}
To finish the proof of Theorem  \ref{cornerthm} it remains to show
that the function $d$ is lower semicontinuous. Assume for the moment
that there exists a sequence
 of points $x_k$ converging to $x$ so that
$d(x_k)>d(x)$.  Since $\varphi$ is continuous, we have  the
convergence $\varphi_1(x_k)=(r^k_1,\ldots ,r^k_n,q^k)\to
\varphi_1(x)=(r_1,\ldots ,r_n,q).$ If for a given  coordinate index
$j$ the coordinate $r^k_j$ vanishes for  all but finitely many $k$,
then $r_j=0$, and if $r^k_j>0$ for all but finitely many $k$, then
$r_j\geq 0$. Hence $d(x_k)\leq d(x)$ contradicting  our assumption.
The proof of  Theorem \ref{cornerthm} is complete.
\end{proof}

\begin{defn} The closure of a connected
component of the set $X(1)=\{x\in X\vert \, d(x)=1\}$ is called a
{\bf face} of the M-polyfold $X$.\index{face}
\end{defn}
Around every point $x_0\in X$ there exists an open neighborhood
$U=U(x_0)$ so that every $x\in U$ belongs to precisely $d(x)$ many
faces of $U$. This is easily verified. Globally it is always true
that $x\in X$ belongs to at most $d(x)$ many faces and the strict
inequality is possible.
\begin{defn} The M-polyfold $X$ is called {\bf face structured},
\index{face structured M-polyfold} if every point $x\in X$ belongs to precisely
$d(x)$ many faces.
\end{defn}
The concept is related to some notion occurring in \cite{L}.

If  $X\times Y$ is a product of two M-polyfolds, then one concludes
 from the definition of the product structure the following relation
 between the degeneracy indices
$$
d_{X\times Y}(x,y)=d_X(x)+d_Y(y).
$$

\subsection{Submanifolds}\label{immsub}
There are many different types of distinguished subsets of a
M-polyfold which qualify as some kind of sub-polyfold. We refer the
reader to \cite{HWZ-polyfolds1} for a comprehensive discussion,
where we introduced three different notions of a sub-polyfold .
Among those one can find sub-polyfolds of locally constant finite
dimensions. These occur as solution sets of nonlinear Fredholm
operators. In this paper we only consider the latter and introduce
the notion of a strong submanifold of a M-polyfold. The more general
notion of a submanifold requires some more work and is given in
\cite{HWZ2}. We just note that both type of submanifolds inherit
from the ambient M-polyfold the structure of a smooth manifold. The
strong submanifolds however lie in a better way in the M-polyfold.

We consider two sc-smooth splicings
\begin{equation*}
{\mathcal S}=(\pi, E,V)\quad \text{and}\quad {\mathcal T}=(\rho,
F,V)
\end{equation*}
 having projections $\pi_v$ and ${\rho}_v$ parametrized by the
same open subset $V$ of a partial quadrant. We define their {\bf
Whitney sum}\index{Whitney sum of splicing} to be the sc-smooth
splicing
$$
{\mathcal S}\oplus{\mathcal T}=(\pi\oplus\rho,E\oplus F,V)
$$
defined by the family of projections
\begin{equation}\label{sect7.1eq2}
(\pi\oplus\rho)_v(h\oplus k)= (\pi_v(h),\rho_v(k)),\quad v\in V.
\end{equation}
One verifies readily that the splicing core $K^{{\mathcal
S}\oplus{\mathcal T}}$ is the fibered sum over $V$ of the splicing
cores $K^{\mathcal S}$ and $K^{{\mathcal T}}$,
\begin{equation}\label{sect7.1eq3}
\begin{gathered}
K^{{\mathcal S}\oplus{\mathcal T}}=K^{\mathcal S}\oplus_{V}K^{{\mathcal T}}=\\
\{(v, h,k))\in V\oplus E\oplus F\vert \, \text{$\pi_v(h)=h$ and
$\rho_v(k)=k$}\}.
\end{gathered}
\end{equation}

\begin{defn}\label{sect7.1def4} The sc-smooth map $f:X\to Y$ between two
M-polyfolds is called a {\bf fred-submersion}\index{fred-submersion}
if
 at every point $x_0\in X$  resp. $f(x_0)\in Y$
there exists a chart $(U, \varphi, {\mathcal T}\oplus \wh{\mathcal
T})$  resp. $(W, \psi, {\mathcal T})$ satisfying  $f(U)\subset W$
and
$$\psi\circ f\circ \varphi^{-1} (v, e', e'')=(v, e')$$
and, moreover, the splicing $\wh{\mathcal T}=(\wh{\rho}, E,V)$ has
the special property that the projections $\wh{\rho}_v$ do not
depend on $v$ and project onto a finite dimensional subspace of $E$.
\end{defn}
Instead of $\wh{\mathcal T}=(\wh{\rho}, E,V)$ we just may take the
splicing $(Id,{\mathbb R}^n,V)$ where $n$ is the dimension of the
image of the projection $\pi$ and $Id$ stands for the constant
family $v\rightarrow Id$. Hence we may assume that in the  Whitney
sum ${\mathcal T}\oplus \wh{\mathcal T}$ the latter summand has the
special form and we will indicate that by writing
 $$
{\mathcal T}\oplus {\mathbb R}^n.
 $$
 The following result will be used quite often.
\begin{prop}\label{fredfred}
If $f:X\to Y$ and $g:Y\to Z$ are fred-submersions,  then the
composition $g\circ f:X\rightarrow Z$ is again a  fred-submersion.
\end{prop}
\begin{proof}
Let $y_0=f(x_0)$ and $z_0=g(y_0)$. We find special charts $\phi$ and
$\psi$ around $x_0$ and $y_0$, respectively, so that
$$
\psi\circ f\circ \phi^{-1}(v,e,e')=(v,e).
$$
Similarly, we find special charts $\alpha$ and $\beta$ so that
$$
\alpha\circ g\circ \beta^{-1}(w,h,h')=(w,h).
$$
Define the inverse of a chart $\gamma$ around $x_0$ by
$$
\gamma^{-1}(w,h,h',e')=\phi^{-1}(\psi\circ\beta^{-1}(w,h,h'),e').
$$
Then we compute
\begin{equation*}
\begin{split}
& \alpha\circ (g\circ f)\circ \gamma^{-1}(w,h,h',e')\\
&= \alpha\circ g\circ f\circ \phi^{-1}(\psi\circ\beta^{-1}(w,h,h'),e')\\
&=\alpha\circ g\circ \psi^{-1}\circ (\psi\circ f\circ \phi^{-1})(\psi\circ\beta^{-1}(w,h,h'),e')\\
&=\alpha\circ g \circ \psi^{-1}\circ (\psi\circ\beta^{-1}(w,h,h'))\\
&=\alpha\circ g\circ\beta^{-1}(w,h,h')\\
&=(w,h).
\end{split}
\end{equation*}
The splicings used for the  charts  involved  are of the form
${\mathcal S}$ and ${\mathcal S}\oplus ({\mathbb R}^n\oplus {\mathbb
R}^k)$. This completes the proof.
\end{proof}
The preimages of smooth points under a fred-submersion carry in a
natural way the structure of smooth manifolds.
\begin{prop}\label{sect7.1prop7.6}
If $f:X\to Y$ is a  fred-submersion  between two M-polyfolds, then
the preimage of a smooth point $y\in Y$,
$$f^{-1}(y)\subset X,$$
carries in a natural way the structure of  a finite dimensional
smooth manifold.
\end{prop}
\begin{proof}
We can define local charts induced from the charts of $X$
(exhibiting $f$ as a fred-submersion). They are defined on open
subsets in ${\mathbb R}^n$. Here $n$ is locally constant, i.e.~only
depends on the connected components of $X$. The transition maps are
sc-smooth and consequently smooth in the classical sense. In other
words, there is a natural system of charts which define the
structure of a smooth manifold on $f^{-1}(y)$.
\end{proof}

The above discussion prompts the following useful concept.

\begin{defn}\label{sect7.1def7.9}
A subset $N\subset X$ of an M-polyfold $X$ is called a {\bf strong
finite-dimensional submanifold} \index{finite-dimensional
submanifold} of $X$ if the following statements hold true.
\begin{itemize}
\item[(i)]\: $N\subset X_{\infty}$.
\item[(ii)]\: For every point $m\in N$ there exists
an open neighborhood $U\subset X$
of $m$, and an M-polyfold \  $Y$, and a surjective fred-submersion
$f:U\to Y$ satisfying
$$f^{-1}(f(m))=N\cap U.$$
\end{itemize}
\end{defn}

The definition of a finite-dimensional submanifold of a M-polyfold
will be given in \cite{HWZ2}.

\section{ M-Polyfold Bundles}\label{chap4}
In this section we continue with the conceptual framework. First we
describe the local models for strong M-polyfold bundles and smooth
maps between them. Then we introduce the notion of a strong
M-polyfold bundle.

\subsection{Local Strong M-Polyfold Bundles}\label{section4.1}
In this subsection we shall introduce the local models for strong
bundles over M-polyfolds. For this we need a generalization of the
notion of splicing where the splicing projection is parameterized by
an open subset of a splicing core. We begin by introducing these
more general splicing definitions.

\begin{defn}\label{vnewdefn4.3}
A  {\bf general sc-smooth splicing}  is a triple
$${\mathcal R}=(\rho,F,(O,{\mathcal S})),$$
where $(O, {\mathcal S})$ is a local M-polyfold model associated
with the sc-smooth splicing ${\mathcal S} =(\pi, E, V)$  and $O$ is
an open subset of the splicing core $K^{{\mathcal S}}=\{(v, e)\in V\oplus E\vert \, \pi_v (e)=e\}$. The space
$F$ is an sc-smooth Banach space and the mapping
\begin{align*}
\rho&: O\oplus F\to F\\
&((v, e), u)\mapsto \rho (v, e, u)
\end{align*}
is sc-smooth. Finally, for fixed $(v, e)\in O$, the mapping
$$\rho_{(v, e)}=\rho (v, e, \cdot ):F\to F$$
is a projection in ${\mathcal L}(F)$. Sc-smoothness of $\rho$, of
course, means that the map
$$
(v,e,u)\rightarrow \rho(v,\pi_v(e),u)
$$
which is defined on an open subset $\wh{O}$ of a partial quadrant
in a sc-Banach space, is sc-smooth.
\end{defn}

The novelty of this definition consists in the requirement, that the
family of projections is parameterized by elements of an open subset
of a splicing core. Iterating this procedure we obtain splicings
parameterized by open sets in splicing cores of generalized
splicings. Continuing this way we arrive at a hierarchy of splicings
of the following types,
\begin{eqnarray*}
&(v,e)\rightarrow (v,\pi_v(e))&\\
&(v,e,u)\rightarrow (v,\pi_v(e),\rho_{(v,\pi_v(e))}(u))&\\
&(v,e,u,w)\rightarrow
(v,\pi_v(e),\rho_{(v,\pi_v(e))}(u),\sigma_{(v,\pi_v(e),\rho_{(v,\pi_v(e))}(u))}(w))&
\end{eqnarray*}
and so on. Hence there are {\bf splicings of type $\boldsymbol{0}$},
which are the original ones, then there are splicings of {\bf type
$\boldsymbol{1}$}, which are the generalized splicings introduced
above, and so on. A type-$k$ splicing can also be viewed as a type
$\ell$-splicing for every $\ell\geq k$. The notion of
$\ssc^r$-smoothness generalizes to these more general splicings.
Using open sets of splicing core of splicings of type $k$ as local
models we can construct {\bf M-polyfolds of type $k$} the same way
we did in Definition \ref{DEF38} for the original M-polyfolds, which
now become M-polyfolds of type $0$. In this paper we shall only meet
M-polyfolds of type $0$ and of type $1$.

The definition of the {\bf tangent of a general sc-smooth splicing}
${\mathcal R}=(\rho,F,(O,{\mathcal S}))$ is defined, quite similarly
as in the case of a splicing, by
$$
T{\mathcal R}=(T\rho,TF,(TO,T{\mathcal S})),
$$
which  is again a general sc-smooth splicing. The map
$T\rho:TO\oplus TF\rightarrow TF$ is a family of projections acting
on $TF$ and parameterized by the tangent $TO$ of $O$. It is defined
by
$$
T\rho(w,\delta w,u, \delta u)=(\rho(w,u),D\rho(w,u)(\delta w,\delta
u)).
$$
Here $w=(v,e)\in O_1\subset V_1\oplus E_1$ and $\delta w\in W\oplus
E$ so that $(w,\delta w)\in TO$ and $(u,\delta u)\in F^1\oplus
F=TF$. Keeping $(w,\delta w)\in TO$ fixed,  the map
$$
T\rho_{(w,\delta w)}:TF\rightarrow TF
$$
is a projection in ${\mathcal L}(F_1\oplus F)$.

Next we introduce the notion of a strong bundle splicing.
\begin{defn}\label{vnewdefn4.3.2}
A {\bf strong bundle splicing} is a general sc-smooth splicing
$${\mathcal R}=(\rho, F, (O, {\mathcal S}))$$
having the following additional property. If $(v,e)\in O_m$ and
$u\in F_{m+1}$, then $\rho ((v,e), u)\in F_{m+1}$ and the  newly
defined triple
$${\mathcal R}^1=(\rho, F^1, (O, {\mathcal S}))$$
is also a general sc-smooth splicing. If we view the strong bundle
splicing ${\mathcal R}$ only as a general smooth splicing we denote
it by ${\mathcal R}^0$.
\end{defn}
Let us note that the complementary splicing ${\mathcal R}^c$ is a
strong bundle splicing as well. From the above definition we
conclude, in particular, that a strong bundle splicing ${\mathcal
R}$ gives rise to two general sc-smooth splicings, namely ${\mathcal
R}^0$ and ${\mathcal R}^1$.

There is a non-symmetric product $E\triangleleft F$ of two sc-Banach
spaces $E$ and $F$. This product is the Banach space $E\oplus F$
equipped, however, with the bi-filtration defined by
$$
(E\triangleleft F)_{m,k}=E_m\oplus F_k
$$
for pairs $(m,k)$  satisfying  $m\geq 0$ and $0\leq k\leq m+1$. For a subset
$U\subset E$ we can define $U\triangleleft F$ in the obvious way.

 The {\bf splicing core}
$K^{\mathcal R}$ of the strong bundle splicing ${\mathcal R}=(\rho,
F, (O, {\mathcal S}))$ is the set
$$
K^{\mathcal R}=\{(w, u)\in O\oplus F\vert\, \rho (w, u)=u\}.
$$
Since ${\mathcal R}$ gives us two general splicings ${\mathcal R}^0$
and ${\mathcal R}^1$ we have a well-defined bi-filtration on
$K^{\mathcal R}$ by pairs $(m,k)$ satisfying $0\leq k\leq m+1$ so
that $K^{\mathcal R}$ can be viewed as a subset of $(V\oplus
E)\triangleleft F$ equipped with the induced  bi-filtration. More
precisely,
$$K^{\mathcal R}_{m, k}=\{(w, u)\in K^{\mathcal R}\vert\,
 w\in O_m,\, u\in F_k\}$$
where $m\geq 0$ and $0\leq k\leq m+1$.  The bundle
$$
K^{\mathcal R}\rightarrow O
$$
defined by means of the strong bundle splicing ${\mathcal R}$ is
called a {\bf local strong  bundle}. It will serve as our local
model of the strong M-polyfold bundles introduced in the next
subsection.
\begin{rem}\label{specialcase}
There is a special case which we already shortly mentioned before.
Assume the strong bundle splicing ${\mathcal R}$ has the special
form
$$
{\mathcal R}= (Id,F,(O,{\mathcal S})),
$$
where ${\mathcal S}=(Id,E,V)$ and $O$ is a relatively open subset of
$V\oplus E$. In this case the splicing core is the product
$$
K^{\mathcal R}=O\triangleleft F
$$
and we can view $O$ as a local model for an sc-manifold and the
product $O\triangleleft F$ as a  model for a local strong sc-bundle
with base $O$.
\end{rem}

 Associated with the strong bundle splicing  ${\mathcal R}$ we have the
splicing cores  $K^{{\mathcal R}^0}$ and $K^{{\mathcal R}^1}$, which
we denote by $K^{\mathcal R}(0)$ and $K^{\mathcal R}(1)$,
respectively. They are equipped with the filtrations
$$
K^{\mathcal R}(0)_m= K^{\mathcal R}_{m,m}\ \hbox{and}\ \
K^{\mathcal R}(1)_m=K^{\mathcal R}_{m,m+1}.
$$
The natural projection $K^{\mathcal R}\rightarrow O:
(w,u)\rightarrow u$ is sc-smooth in the sense that the two
projections
$$
K^{\mathcal R}(i)\rightarrow O
$$
are sc-smooth for $i=0,1$.

  We can define  the tangent $T{\mathcal R}$
of the strong bundle splicing
$$
{\mathcal R}=(\rho, F, (O, {\mathcal
S}))
$$
 as follows. First we consider the underlying strong bundle
splicing ${\mathcal R}^0$ and take the associated tangent splicing
$T{\mathcal R}^0$,
$$
T{\mathcal R}^0=(T\rho, TF, (TO, T{\mathcal S})).
$$
Since we also have the splicing ${\mathcal R}^1$, we also can take
its tangent $T{\mathcal R}^1$ given by
$$
T{\mathcal R}^1=(T\rho, T(F^1), (TO, T{\mathcal S})).
$$
From $T(F^1)=(TF)^1$ we conclude  that
$$
T{\mathcal R}=(T\rho,TF,(TO,T{\mathcal S}))
$$
is again a strong bundle splicing in the sense of Definition
\ref{vnewdefn4.3.2}. Its {\bf splicing core} $K^{T{\mathcal R}}$ is,
as usual,  defined by
$$
K^{T{\mathcal R}}=\{(w, \delta w, u, \delta u)\in TO\oplus TF\ \vert
\, T\rho (w, \delta w, u, \delta u)=(u,\delta u)\}.
$$
More explicitly, the elements of $K^{T{\mathcal R}}$  are restricted
by the following equations for $w=(v,e)\in O_1\oplus E_1$ and
$\delta w=(\delta v,\delta e)\in W\oplus E$ so that $(w,\delta w)\in
TO$ and for $(u,\delta u)\in TF$,
\begin{equation*}
\begin{gathered}
\pi (v, e)=e\\
\rho (w, u)=u\\
\delta e=\pi (v, \delta e)+D_{v}\pi (v, e)\delta v\\
\delta u=\rho (w, \delta u)+D_w\rho (w, u)\delta w.
\end{gathered}
\end{equation*}
Let us observe that for $i=0,1$ the following relationships hold for
the underlying general sc-smooth splicings
$$
(TK^{\mathcal R})(i) = K^{T{\mathcal R}}(i)  = K^{T({\mathcal
R}^i)}=T(K^{\mathcal R}(i)).
$$

 Next we shall define the concept of a strong bundle map of class
 $\ssc^k_\triangleleft$
between splicing cores of strong bundle splicings. Recall the
Definition \ref{def2.7} for the $\ssc^1$-class of mappings between
open subsets of splicing cores.

\begin{defn}\label{vnewdefn4.3.3}({\bf Strong bundle maps})
If ${\mathcal R}=(\rho, F, (O, {\mathcal S}))$ and ${\mathcal
R}'=(\rho', F', (O', {\mathcal S}'))$ are two strong bundle
splicings we denote the associated splicing cores by $K=K^{\mathcal
R}\subset O\oplus F$ and $K'=K^{{\mathcal R}'}\subset O'\oplus F'$.
Consider a map
$$f:K\to K'$$
of the form
$$f(w, u)=(\varphi (w), \Phi (w, u)),$$
where $\varphi:O\to O'$ and where $\Phi :O\oplus F\to F'$. Then
\begin{itemize}
\item[$\bullet$]\: The map $f$ is  a {\bf strong bundle map of
class} {\bf $\ssc^{\mathbf 0}_{\mathbb \trl}$}, or simply an
$\ssc^0_\triangleleft$-map if it induces $\ssc^0$-maps $K(i)\to
K'(i)$ for $i=0$ and $i=1$.
\item[$\bullet$]\: The map $f$ is a  {\bf strong bundle map} of
{\bf class ${\boldsymbol \ssc}_{\boldsymbol{\trl}}^{\mathbf 1}$} if
it is of class $\ssc^0_{\trl}$ and if  it induces $\ssc^1$-maps
$K(i)\to K'(i)$ for $i=0$ and $i=1$.
\end{itemize}
Observe that $K(i)$ and $K'(i)$ are type-$1$ M-polyfolds.
\end{defn}
In many cases we require $\Phi$ to be linear in $u$. In particular,
this is the case when $\Phi$ occurs as an isomorphism between local
strong M-polyfold bundles.

Next we  consider  maps $f:K\rightarrow K'$ between splicing cores of
strong bundle splicings of the form as in Definition
\ref{vnewdefn4.3.3}. In order to define maps $f:K\rightarrow K'$ of
class $\ssc_\triangleleft^2$ we proceed as in the sc-case. Assuming
that $f$ is of class $\ssc_\triangleleft^1$ we consider it first as
an $\ssc^1$-map
$$
f:K(0)\rightarrow K'(0)
$$
between splicing cores. Its tangent map $Tf$ is described by the
formula
$$
Tf(w,\delta w,u,\delta u)=(\varphi(w),D\varphi(w)\delta
w,\Phi(w,u),D\Phi(u,w)(\delta w,\delta u))
$$
Since $f$ is also an $\ssc^1$-map
$$
f:K(1)\rightarrow K'(1)
$$
the tangent formula above defines two maps
$$
Tf:TK^{\mathcal R}(i)=K^{T{\mathcal R}^i}\rightarrow TK^{{\mathcal
R}'}(i)=K^{T{{\mathcal R}'}^i}
$$
for $i=0,1$ which are $\ssc^0$-continuous. Therefore they define a
map
$$
Tf:TK^{\mathcal R}\rightarrow TK^{{\mathcal R}'}
$$
between splicing cores of strong bundle splicings which is of class
$\ssc^0_\triangleleft$. It is called the {\bf tangent map} of the
$\ssc^1_\triangleleft$-map $f$.

 If this tangent map $Tf$ is of class
$\ssc^1_{\trl}$ as defined above, then the map $f:K\to K'$ is called
of {\bf class $\ssc^{\mathbf 2}_{\mathbf \trl}$}. Proceeding
inductively as in the sc-case one defines the mappings $f:K\to K'$
of class $\ssc^k_{\trl}$ for $k\geq 1$ and also the
$\ssc_{\trl}$-smooth mappings. Let us finally note that the chain
rule also holds for strong bundle maps.
\begin{thm}({\bf Chain rule for strong bundle maps})
Let $f:K|O\rightarrow K'$ and $g:K'|O'\rightarrow K''$ be two strong
bundle maps of class $\ssc^1_\triangleleft$ between local strong
bundles so that the image of $f$ is contained in the domain of $g$.
Then the composition $g\circ f$ is also a strong bundle map of class
$\ssc^1_\triangleleft$ and the tangent maps satisfy
$$
T(g\circ f) =Tg\circ Tf.
$$
\end{thm}
A  note of caution! As before the order of terms in the tangent map
$Tf$ of a $\ssc^1_\triangleleft$-map is different from their order
in the classical theory.

Associated with  the strong bundle splicing ${\mathcal R}$  we have
the local strong  bundle $p:K\rightarrow O$. A sc-smooth section of
the bundle $p$ is just an sc-smooth section of the underlying bundle
$K(0)\rightarrow O$. The vector space of sc-smooth sections is
denoted by $\Gamma(p)$. In addition, there is a different class  of
sections called $\ssc^+$-sections. An $\ssc$-smooth section $f$ is
called a {\bf $\ssc^+$-section}  if it defines an sc-smooth section
of the bundle $K(1)\rightarrow O$. We denote the collection of
$\ssc^+$-sections by $\Gamma^+(p)$.

\subsection{Fillability and Fillers}\label{section4.2}
Considering the local strong bundle $p:K^{\mathcal R}\to O$
associated with a strong bundle splicing ${\mathcal R}$, we
investigate the coherence in the jumps of the space dimensions in
the base and the fibers.

We start with a strong bundle splicing ${\mathcal
R}=(\rho,F,(O,{\mathcal S}))$. The splicing ${\mathcal S}$ is  the triple
 $(\pi,E,V)$ in which  $V$ is an open subset  of  a partial
quadrant $C$ contained in the sc-Banach space $W$. The set $O$ is an open neighborhood of  the origin  in
the splicing core $K^{\mathcal S}=\{(v, e)\in V\oplus E\vert \, \pi_v (e)=e\}$.
 If
$$
s:O\rightarrow V
$$
is  the sc-smooth map defined by $s(v, e)=v$, we shall abbreviate by
$s^*\pi:O\times E\to E$ the composition $s^*\pi (v, e, u)=\pi(s(v,
e), u)=\pi (v, u)$ and introduce the general sc-smooth splicing
 $s^\ast{\mathcal S}^c$, having the splicing parameter set $O$, by
 $$s^\ast{\mathcal S}^c=(1-s^*\pi, E, O).$$
 Its splicing core is the set
 \begin{equation*}
 \begin{split}
 K^{s^\ast{\mathcal S}^c}&=\{(v, e, u)\in O\oplus E\vert \, \pi_{s(v, e)}(u)=0\}\\
 &=\{(v, e, u)\in V\oplus E\oplus  E\vert \, \text{$(v, e)\in O, \, \pi_v(e)=e$ and $(1-\pi_v)(u)=u$}\}.
 \end{split}
 \end{equation*}
 In view of the splitting $E=\pi_v (E)\oplus (1-\pi_v)(E)$, the splicing core $ K^{s^\ast{\mathcal S}^c}$ can be naturally identified with the following open subset  $\wh{O}$ of $V\oplus E$,
 $$
\wh{O}=\{(v,e)\in V\oplus E\ |\ (v,\pi_v(e))\in O\}.
$$
We have the natural projection
$$
a:K^{s^\ast{\mathcal S}^c}\rightarrow O, \qquad (v,w)\rightarrow
(v,\pi_v(w))
$$
and we can view  $a:K^{s^\ast{\mathcal S}^c}\rightarrow O$  as a  bundle (of course, not as a strong bundle). The fiber $a^{-1}(v, e)$ over the  point $(v,e)\in O$ is the Banach
space
\begin{equation*}
\begin{split}
a^{-1}(v,e)&=\{(v,e+w)\vert \, w\in E, \pi_v(w)=0\}\\
&=\{(v,e)\}\times \ker(\pi_v).
\end{split}
\end{equation*}
The strong bundle splicing ${\mathcal R}$ comes together with its
complementary strong bundle splicing ${\mathcal R}^c=(1-\rho, F, (O,
{\mathcal S}))$ giving rise to   the local strong bundle
$b:K^{{\mathcal R}^c}\rightarrow O$. We are interested  only in the
underlying bundle
$$
b:K^{{\mathcal R}^c}(0)\rightarrow O
$$
associated with the general sc-smooth splicing ${({\mathcal R}^c)}^0$.

The following concept of a filler turns out to be very useful in the applications.
\begin{defn}\label{ndef4.5}
Let ${\mathcal R}$ be a strong bundle splicing and ${\mathcal R}^c$
the associated complementary strong bundle splicing.
Consider the two bundles over $O$
$$
a:K^{s^\ast {\mathcal S}^c}\rightarrow O\ \ \hbox{and}\ \
b:K^{{\mathcal R}^c}(0)\rightarrow O.
$$
Then a {\bf filler}  for ${\mathcal R}$  is an  sc-diffeomorphism
$$
f^c: K^{s^\ast{\mathcal S}^c}\rightarrow K^{{\mathcal R}^c}(0)
$$
between the complementary bundle pairs,
which is linear in the fibers and covers the identity map
 $O\rightarrow O$. (It is, in particular, a bundle  isomorphism).
\end{defn}
\begin{defn}\label{ndefn4.6}
The strong bundle splicing ${\mathcal R}$ is {\bf fillable}
if there
exists a filler for ${\mathcal R}$.
\end{defn}
Being fillable is a property of the strong bundle splicing
${\mathcal R}$.

A filler $f^c:K^{s^\ast{\mathcal S}^c}\rightarrow
K^{{\mathcal R}^c}(0)$ has the  the form
$$
f^c: (v,e,u)\rightarrow (v,e,{\bf f}^c(v,e,u))
$$
 where $(v,e)\in O$ and where
$u\in E$  satisfies $\pi_v(u)=0$.
The principal part  ${\bf f}^c (v, e, u)\in F$
satisfies
 $\rho_{(v,e)}({\bf f}^c(v,e,u))=0$.
 In view of the identity $e=\pi_v (e)+(1-\pi_v )(e)$ in $E$,
 the principal part ${\bf f}^c$ can be viewed  as an sc-smooth  map
 $$
 \wh{O}\rightarrow F, \qquad (v,e)\rightarrow {\bf f}^c(v,e)
$$
satisfying $\rho_{(v,\pi_v(e))}({\bf f}^c(v,e))=0$.

\subsection{Strong M-Polyfold Bundles}\label{sect4.4}
In order to introduce strong M-polyfold bundles we consider a
surjective sc-smooth map
$$p:Y\to X$$
between two M-polyfolds. The M-polyfold $Y$ is of type-$1$ and $X$
of type-$0$. We assume in addition that for every $x\in X$, the
preimage $p^{-1}(x)=Y_x$, called the fiber over $x$,
 carries the structure
of a Banach space.
\begin{defn}\label{def4.7}
Let $p:Y\rightarrow X$ be as just described. A  {\bf strong  M-polyfold
bundle chart} for the bundle $p:Y\rightarrow X$ is a triple
$(U,\Phi, (K^{\mathcal R},{\mathcal R}))$. Here $U\subset X$ is an
open set and ${\mathcal R}=(\rho , F, (O, {\mathcal S}))$  a strong
bundle splicing with the local model $(O, {\mathcal S})$  of the
polyfold $X$. The map $\Phi$ is an  sc-diffeomorphism
$$
\Phi:p^{-1}(U)\rightarrow K^{\mathcal R}
$$
which is linear on the fibers and  covers   the   sc-diffeomorphism
$$
\varphi:U\rightarrow O
$$
so that $pr_1\circ\Phi=\varphi\circ p$. Moreover,  the map $\Phi$
resp. $\varphi$ are smoothly compatible with the M-polyfold
structures on $Y$ and  $X$, respectively.

 \begin{equation*}
 \begin{CD}
 p^{-1}(U)@>\Phi >> K^{\mathcal R}\\
           @VpVV       @VV\text{pr}_1V       \\
 U@>\varphi >> O.
  \end{CD}
 \end{equation*}
 \mbox{}\\[1ex]
\end{defn}

Recall that ${\mathcal S}=(\pi , E, V)$ is an sc-smooth splicing
where $V$ is an open subset of a partial quadrant in an sc-smooth
Banach space $W$. The set $O$ is an open subset of the splicing core
$K^{\mathcal S}=\{w=(v, e)\in V\oplus E\vert \, \pi (v, e)=e\}$
while the splicing core $K^{\mathcal R}$ over the base $O$ is
defined by $K^{\mathcal R}=\{(w, u)\in O\oplus F\vert \, \rho (w,
u)=u\}$ where $F$ is an sc-smooth Banach space.

\begin{defn}
Two M-polybundle charts $(\Phi, \varphi )$ and $(\Psi, \psi)$  are
called {\bf $\mathbf{sc}_{\trl}$-compatible},  if
 the transition map
 $$\Psi \circ \Phi^{-1}:K^{\mathcal R}\vert \varphi (U\cap U')\to
 K^{{\mathcal R}'}\vert \psi (U\cap U')$$
 between their splicing cores $K^{\mathcal R}$ and $K^{{\mathcal R}'}$
 is an $\ssc_{\trl}$-smooth strong bundle map.

 An {\bf M-polybundle atlas}  consists of a
family of M-polybundle charts $(U, \Phi, (K^{\mathcal R}, {\mathcal
R}))$ so that  the underlying open sets $U$ cover $X$ and so that
the transition maps are sc$_\trl$-smooth strong bundle maps. A
maximal smooth atlas of M-polybundle charts is called an {\bf
M-polybundle structure} and the bundle
$$p:Y\rightarrow X$$
is called a {\bf strong M-polyfold bundle.}

\end{defn}

\begin{defn}\label{ndefn4.8}
A strong $M$-polyfold bundle $p:Y\to X$ is called
 {\bf fillable} if around every point $q\in X$,
 there exists a compatible strong $M$-polyfold
 bundle chart $(U, \Phi, (K^{{\mathcal R}}, {\mathcal R}))$
 whose strong bundle  splicing ${\mathcal R}$ is fillable.
\end{defn}

As it turns out all strong bundles occurring in the applications we
have in mind have this property.

Note that the tangent bundle $TX\rightarrow X^1$ in general is not a
strong M-polyfold bundle.

 Given the strong polyfold bundle $p:Z\to X$ having the base $X$ and
 the sc-smooth map $f:Y\to X$ between M-polyfolds one defines the
 (algebraic) pullback (or induced) bundle
$$p': f^*Z\to Y$$
having the base $Y$ as follows. One  takes  the set $f^*Z=\{ (y,
z)\in Y\times Z\vert \, p(z)=f(y)\}$ and the two projection maps $p'
(y, z)=y$ and $f'(y, z)=z$, so that the diagram

\begin{equation*}
 \begin{CD}
f^*Z@>f' >>Z\\
           @Vp'VV       @VVpV      \\
Y@>f> > X
 \end{CD}
 \end{equation*}
 \mbox{}\\[1ex]
 commutes.

\begin{prop}\label{newpropA}
If $p:Z\to X$ is a strong M-polyfold bundle and $f:Y\rightarrow X$
is an sc-smooth map between M-polyfolds, then the pullback bundle
$p':f^*Z\to Y$ carries a natural induced structure of a strong
M-polyfold bundle whose base is the M-polyfold $Y$.
\end{prop}
\begin{proof}
Choose a point $(y_0, z_0)\in f^\ast Z$ so that
$f(y_0)=p(z_0)=x_0\in X$. Take a strong $M$-polyfold bundle chart for $p:Z\to
X$ denoted by $(U, \Phi, (K^{\mathcal R}, {\mathcal R}))$,
 \begin{equation*}
 \begin{CD}
 p^{-1}(U)@>\Phi >> K^{\mathcal R}\\
           @VpVV       @VV\text{pr}_1V       \\
 U@>\varphi >> O,
  \end{CD}
 \end{equation*}
 \mbox{}\\[1ex]
so that  the open set $U\subset X$ contains the point $x_0$. The
strong bundle splicing ${\mathcal R}=(\rho , F, (O, {\mathcal S}))$
is associated with the local model $(O, {\mathcal S})$ of the
polyfold $X$, where $O$ is an open subset of the splicing core
$K^{\mathcal S}$ of the splicing ${\mathcal S}=(\pi , E, V)$. Take
now an M-polyfold chart $\psi :U'\to O'$ around the given point
$y_0\in U'\in Y$, associated  with the  local model $(O', {\mathcal
S}')$ of the M-polyfold $Y$, where $O'\subset K^{{\mathcal S}'}$ is
an open subset of the splicing core of the sc-smooth splicing
${\mathcal S}'=(\pi', E', V')$. Choose $U'\subset Y$ so small that
$f(U')\subset U$. Define the strong bundle splicing ${\mathcal
R}'=(\rho', F, (O', {\mathcal S}'))$ by means of the sc-smooth map
$\rho':O'\oplus F\to F$ given as
$$\rho' (v', u):=\rho (\varphi \circ f\circ \psi^{-1}(v'), u).$$
The strong  M-polyfold bundle chart for the bundle $p':f^*Z\to Y$ is
now the triple $(U',  \Psi , (K^{{\mathcal R}'}, {\mathcal R}'))$
with the homeomorphism
$$\Psi:f^*Z\vert U'\to K^{{\mathcal R}'}\vert O'$$
defined as
$$\Psi (p'^{-1}(y))=\Phi (p^{-1}(f(y))),$$
for all $y\in U'\subset Y$. We have the diagram
 \begin{equation*}
 \begin{CD}
 p'^{-1}(U')@>\Psi >> K^{{\mathcal R}'}\\
           @Vp'VV       @VV\text{pr}_1V       \\
 U'@>\varphi >> O'.
  \end{CD}
 \end{equation*}
 \mbox{}\\[1ex]
One can verify that the transition maps between two such strong
M-polyfold bundle charts are $\ssc_{\trl}$-smooth. This finishes the
proof of Proposition \ref{newpropA}
\end{proof}

\subsection{Sections and Linearizations}
Assume that $p:Y\rightarrow X$ is a strong M-polyfold bundle over
the M-polyfold $X$. We denote the space of sc-smooth sections by
$\Gamma(p)$. In addition there is the distinguished space of
$\ssc^+$-sections which we denote by $\Gamma^+(p)$.

The section $f$ is an {\bf $\ssc^+$-section} if its local
representations in the strong $M$-polyfold bundle charts are
$\ssc^+$-sections as defined at the end of section \ref{section4.1}
above. If $(U, \Phi, (K^{{\mathcal R}}, {\mathcal R}))$ is  such a
strong $M$-polyfold bundle chart, the local representation of the
section $f$ of the bundle $p:Y\to X$ is the section $g$ of the
strong local bundle $K^{{\mathcal R}}\to O$ defined as the
push-forward of $f$ by
$$g(w)=\Phi \circ f\circ \varphi^{-1}(w),$$
where $w\in O$. By definition, the map $\Phi:p^{-1}(U)\to
K^{{\mathcal R}}$ is an sc-diffeomorphsim  which is linear in the
fibers and which covers the sc-diffeomorphism $\varphi:U\to O$ where
$O$ is an open subset of the splicing core $K^{{\mathcal S}}=\{(v,
e)\vert \, \pi_v (e)=e\}$ belonging to the splicing ${\mathcal
S}=(\pi, E, V)$. Associated with the strong bundle splicing
${\mathcal R}= (\rho, F, (O, {\mathcal S}))$ we have the splicing
core $K^{{\mathcal R}}=\{(w, u)\in  O\oplus F\vert \, \rho (w,
u)=u\}$.

Next we choose  a smooth point $q\in X$. Generalizing a trivial classical fact for vector bundles  we can identify naturally the tangent space $T_{0_q}Y$ at the zero element $0_q=\Phi^{-1}(\varphi (q), 0_F)$ with the sc-Banach space $T_qX\oplus Y_q$ where $Y_q=p^{-1}(q)$ is the fiber.  Since $q$ is a smooth point we may assume in the following that $\varphi (q)\in O$ is equal to $0$, so that $\Phi (0_q)=(0, 0)\in K^{{\mathcal R}}\subset O\oplus F$. The identification $T_{0_q}Y
\longleftrightarrow T_qX\oplus Y_q$ corresponds in the local coordinates to the identification
$$(0, \delta w, 0, \delta u)\longleftrightarrow ((0, \delta w), (0, \delta u))$$
of the elements in the tangent space $T_{(0, 0)}K^{{\mathcal
R}}\subset T_0O\oplus T_0F$. We shall denote by $P_q:T_{0_q}Y\simeq
T_qX\oplus Y_q\to Y_q$ the projection.

Given a section $f\in\Gamma(p)$ which vanishes at the  smooth point $q\in X$ we
define, following the classical recipe, the {\bf linearization}
$f'(q)$ by
$$
f'(q):T_qX\rightarrow Y_q:h\rightarrow P_q\circ Tf(q)h.
$$
As in the case of vector bundles there is generally not an intrinsic
notion of a linearization of a section at an arbitrary point $q$ at
which $f(q)$ does not vanish. In our case with $Y\rightarrow X$
being a strong bundle we have, however, some additional structure.
This will allow us to define a linearization at an arbitrary smooth
point which is unique up to a linear $\ssc^+$-operator.

In order to see this, we consider the sc-smooth section $f\in \Gamma (p)$ and look at the smooth point $q\in X$. Its image $y=f(q)$ is a smooth point in $Y$ and we claim that there exists an $\ssc^+$-section $s$ defined near $q$ and satisfying $s(q)=f(q)$.

Indeed, if the coordinate representation of the section $f$ is given
by $g(w)=(w, {\bf g}(w))$ and if $q$ corresponds to $w_0$, then
${\bf g}(w_0)$ is a smooth point in the sc-Banach space $F$
satisfying $\rho (w_0, {\bf g}(w_0))={\bf g}(w_0)$. Now define in
the local coordinates the section $s$ by $s(w)=(w, \rho (w, {\bf
g}(w_0))$. It satisfies $s(w_0)={\bf  g}(w_0)$ and is indeed an
$\ssc^+$-section because the projections $\rho$ belong to a strong
bundle splicing as defined in Definition \ref{vnewdefn4.3.2}.

Now take any $\ssc^+$-section $s$ of the bundle $Y\to X$ defined
near $q$ and satisfying $s(q)=f(q)$. Then the section $f-s$ is
defined near $q$ and vanishes at $q$. We define the linearization
$f'_s(q)$ by
$$
f'_s(q):T_qX\rightarrow Y_q:h\rightarrow P_q\circ T(f-s)(q)h.
$$
Next we investigate to what extend $f'_s(q)$ depends on the choice
of $s$. Assume therefore that $s$ and $t$ are $\ssc^+$-sections
defined near $q$ and satisfying  $s(q)=t(q)=f(q)$. Then, by
definition,
\begin{equation*}
\begin{split}
f'_s(q)&=P_q\circ T(f-s)(q)\\
&=P_q\circ T(f-t +(t-s))(q)\\
&=P_q\circ T(f-t)(q) +P_q\circ T(t-s)(q)\\
&= f'_t(q) + P_q\circ T(t-s)(q).
\end{split}
\end{equation*}
It remains to understand the perturbation term $P_q\circ T(t-s)(q)$.
For this we observe that it suffices to understand $P_q\circ Ts(q)$
for an  $\ssc^+$-section $s$ defined in a neighborhood $U\subset X$
of  $q$  and vanishing at $q$. Since $s$ is an  $\ssc^+$-section of
the bundle $p:Y|U\rightarrow U$ its tangent $Ts$ is an
$\ssc^+$-section of $Tp:T(Y|U)\rightarrow TU$. Hence the composition
$$
P_q\circ Ts(q):T_qX\rightarrow Y_q
$$
which in local coordinates is given by
$$
(0,\delta w)\rightarrow (0,\delta w, 0,Ds(0)\delta w)\rightarrow
(0,Ds(0)\delta w)
$$
is an  $\ssc^+$-operator. \begin{defn}
Let $[f,q]$ be the germ of a section $f$ of the strong bundle
$p:Y\rightarrow X$ around the smooth point $q$. Let $[s]$ be a germ
of a $\ssc^+$-section around $q$ which satisfies $s(q)=f(q)$. Then
the {\bf linearization} of $[f,q]$ with respect to $[s]$ is defined
by
$$
f'_{[s]}(q)=P_q\circ T(f-s)(q).
$$
\end{defn}
The above discussion is now summarized in the following proposition.
\begin{prop}\label{ox}
Let $[f,q]$ be an sc-smooth section germ of the bundle
$p:Y\rightarrow X$ near a smooth point $q$. Then two linearizations
$f'_{[s]}(q)$ and $f'_{[t]}(q)$ differ by a $\ssc^+$-operator. In
particular,  if one linearization is sc-Fredholm so are all others.
\end{prop}
The last statement follows from Proposition \ref{prop1.21}.  This
allows us to introduce  the following definition.
\begin{defn}
An sc-smooth section of the strong M-polyfold bundle $p:Y\rightarrow
X$ is {\bf linearized Fredholm at the smooth point $q$} provided a
linearization at the point $q$ is sc-Fredholm. We say $f$ is {\bf
linearized Fredholm} if this holds at all smooth points $q$.
\end{defn}
If $f$ is linearized Fredholm  and $q$ a smooth point, we define
the index $\ind(f,q)\in {\mathbb Z}$ by
$$
\ind(f,q):= i(f'_{[s]}(q)).
$$
In view of Proposition \ref{ox} this is well-defined. Here $i$
denotes the Fredholm index.

Another consequence of the previous discussion is the following.
\begin{prop}
Assume that $p:Y\rightarrow X$ is a strong M-polyfold bundle and
$f\in\Gamma(p)$ an sc-smooth section which is linearized Fredholm.
Then the section $f+s$ for any $\ssc^+$-section in $\Gamma^+(p)$ is
linearized Fredholm.
\end{prop}
If a strong M-polyfold bundle $p:Y\rightarrow X$ is fillable one can
construct for every  $\ssc$-smooth section $f$ near a smooth point
$q\in X$ a filled section. This will be important for the Fredholm
theory developed in \cite{HWZ2}. To carry out this construction we
assume that $f$ is an $\ssc$-smooth section of the bundle
$p:Y\rightarrow X$ and $q\in X_\infty$. We pick a fillable strong
bundle coordinate
$$
\Phi:Y|U\rightarrow K^{\mathcal R}
$$
defined on an open neighborhood $U$ of $q$ and covering the
 sc-diffeomorphism $\varphi:U\rightarrow O$. Let
$$
f^c:K^{s^\ast{\mathcal S}^c}\rightarrow K^{{\mathcal R}^c}(0)
$$
be a filler for the strong bundle splicing  ${\mathcal R}=(\rho,F,(O,(\pi,E,V)))$.
The principal
part of $f^c$ gives us a  sc-smooth map
$$
{\bf f}^c:\wh{O}\rightarrow F:(v,e)\rightarrow {\bf f}^c(v,e)
$$
satisfying $\rho_{(v,\pi_v(e))}({\bf f}^c(v,e))=0$.  Recall that
$\wh{O}$ stands for the open  subset  of $V\oplus E$ defined by
$\wh{O}=\{(v, e)\in V\oplus E\vert \, (v, \pi_v (e))\in O\}$. The
push-forward $\Phi_\ast f$  is a section of the local strong bundle
$K^{{\mathcal R}}\to O$. Its principal part has a natural extension
to the open set $\wh{O}$ which we denote by ${\bf f}$. It satisfies
$$
\rho_{(v,\pi_v(e))} ({\bf f}(v,e))={\bf f}(v,e).
$$
Finally,  we  introduce the sc-smoth map $\ov{\bf f}:\wh{O}\rightarrow F$ by
$$
\ov{\bf f}(v,e)={\bf f}(v,e)+{\bf f}^c(v,e).
$$
which can be  viewed as the principal part of an  sc-smooth
section $\ov{f}$ of the bundle $\wh{O}\triangleleft F\rightarrow \wh{O}$.
\begin{defn}
The section $\bar{f}$ is called a {\bf filled version} of $f$ near
$q$ and $\ov{\bf f}$ is called its principal part, i.e.
$$
\ov{f}(v,e)=((v,e),\ov{\bf f}(v,e)).
$$
\end{defn}
Let us observe that $\ov{f}(v,e)=0$ if and only if ${\bf f}(v,e)=0$
and ${\bf f}^c(v,e)=0$. Since $f^c$ is a filler we deduce from ${\bf f}^c(v,e)=0$ that
 $(1-\pi_v)(e)=0$  implying $\pi_v(e)=e$  so that $(v,e)\in O$
and $\Phi_\ast f(v,e)=0$. Hence $\varphi^{-1}(v,e)$ is a zero of the section
$f$.  Consequently,  a filled version still describes in local coordinates
precisely the solution set of $f$ over the set $U$.

Assume again that $p:Y\rightarrow X$ is a fillable strong M-polyfold
bundle and $f$ an sc-smooth section. Suppose that $\ov{f}$ is a
filled version representing the section $f\vert U$ as a section of
the bundle $\wh{O}\triangleleft F\rightarrow \wh{O}$. We will prove
the following result.
\begin{prop}
For a smooth point $q\in U$ corresponding to the point $(v,e)\in O$
the linearization $f'_{[s]}(q):T_qX\rightarrow Y_q$ is sc-Fredholm
if and only if the linearization of the filled version
$\ov{f}'_{[t]}(v,e):T_{(v,e)}\wh{O}\rightarrow F$ is
$\ssc$-Fredholm. In this case the  Fredholm indices are the same.
\end{prop}
\begin{proof}
Without loss of generality we may assume that $f(q)=0$ and $s=0$.
Using fillable strong bundle coordinates $\Phi:Y\vert U\rightarrow
K^{\mathcal R}$ covering $\varphi:U\rightarrow O$ we may assume
without loss of generality that $f$ is a sc-smooth section of the bundle
$$
p:K^{\mathcal R}\rightarrow O
$$
satisfying  $0\in O$ and $f(0)=0$. Then $f$ has the form
$$
O\rightarrow K^{\mathcal R}:(v,e)\rightarrow ((v,e),\wh{f}(v,e)).
$$
Using its principal part $\wh{f}$ we can define an $\ssc$-smooth
map
$$
{\bf f}:\wh{O}\rightarrow F:(v,e)\rightarrow \wh{f}(v,\pi_v(e)),
$$
which satisfies $\rho_{(v,\pi_v(e))}({\bf f}(v,e))={\bf f}(v,e)$. The open subset $\wh{O}$ of $V\oplus E$ can be naturally identified with the splicing core $K^{s^\ast{\mathcal S}^c}$, as we have seen in section \ref{section4.2}.  Since the strong bundle splicing ${\mathcal
R}$ is fillable,  we have the existence of a bundle isomorphism
(linear in the fibers)
$$
K^{s^{\ast}{\mathcal S}^c}\rightarrow K^{{\mathcal R}^c}(0)
$$
over the set $O$ which gives us the  $\ssc$-smooth  filler
$$
{\bf f}^c:\wh{O}\rightarrow F
$$
satisfying $\rho_{(v,\pi_v(e))}({\bf f}^c(v,e))=0$. In addition, for
fixed $(v,e)\in O$ the map $r\rightarrow {\bf f}^c(v,e+r)$ is a
linear isomorphism between the Banach spaces $\ker(\pi_v)$ and
$\ker(\rho_{(v,e_v)})$. Finally,  the principal  part of the locally
filled section is the sc-smooth map
$$
\ov{\bf f}:\wh{O}\rightarrow F:(v,e)\rightarrow {\bf f}(v,e)+{\bf
f}^c(v,e).
$$
Since $f(0, 0)=0$, the linearisation  of the section $f$ of the bundle $K^{{\mathcal R}}\to O$ at the point  $(0,0)\in O$,
$$
f'(0,0)):T_{(0,0)}O\rightarrow \ker(1-\rho_{(0,0)})=p^{-1}(0,0)
$$
is  equal to $D{\bf f}(0,0)|T_{(0,0)}O$. We have to compare
it to the linearisation $D{\bar{\bf f}}(0,0):W\oplus E\rightarrow
F$. Here $W$ is the sc-Banach space containing the relatively open
neighborhood $V$ of $0$ in a partial quadrant $C\subset W$.
 According to the splitting $E=\pi_0 (E)\oplus (1-\pi_0)(E)$
  we shall split  $\delta e\in E$ into $\delta e=(\delta a, \delta
  b)$ so that the tangent space at the point
  $(v,e)=(0,0)\in O$ becomes $T_{(0,0)}O=\{(\delta w,\delta e)\in
  W\oplus E\ |\ \pi_0(\delta e)=\delta e\}=\{(\delta w,\delta a)\in
  W\oplus \pi_0(E)\}$
  and compute,
\begin{equation*}
\begin{split}
&D{\bar{\bf f}}(0,0)(\delta w,\delta e)\\
&\phantom{=}=D{\bf f}(0,0)(\delta w,\delta e) +D{\bf f}^c(0,0)(\delta w,\delta
e)\\
&\phantom{=}=D{\bf f}(0,0)(\delta w,\delta a)+D{\bf f}(0,0)(0,\delta b)+D{\bf
f}^c(0,0)(\delta w,\delta a)+D{\bf f}^c(0,0)(0,\delta b)\\
&\phantom{=}=D{\bf f}(0,0)(\delta w,\delta a)+D{\bf
f}^c(0,0)(\delta w,\delta a)+D{\bf f}^c(0,0)(0,\delta b)\\
&\phantom{=}=:f'(0,0)(\delta w,\delta a)+B(\delta w,\delta
a)+C(\delta b).
\end{split}
\end{equation*}
We have concluded from the identity ${\bf f}(v, e)={\bf f}(v, \pi_v
(e))$ that
$$D{\bf f}(0, 0)(0, \delta b)=0.
$$
 In addition, since for fixed $v$ the map ${\bf f}^c (v, \pi_v
(e)+(1-\pi_v )(e))$ is linear in $(1-\pi_v )(e)$, we conclude that
$D{\bf f}^c (0, 0)(0, \delta b)={\bf f}^c (0, \delta b)$. Since
$(0,0)\in O$ is a smooth point, the map
$$
C:\ker(\pi_0)\rightarrow \ker(\rho_{(0,0)}):\delta b\rightarrow {\bf
f }^c(0,\delta b)
$$
is a linear $\ssc$-isomorphism,  by the definition of a filler. The
sc-operator $W\oplus \ker(\pi_0)\rightarrow \ker(\rho_{(0,0)})$
defined by $(\delta w,\delta a)\rightarrow D{\bf f}^c(0,0)(\delta
w,\delta a)$ vanishes.

Hence our total linear $\ssc$-operator $D\ov{{\bf f}}(0, 0)$ has the
matrix  form
$$
\begin{bmatrix}
(\delta w,\delta a)\\
\delta b
\end{bmatrix}
\rightarrow
\begin{bmatrix}
f'(0,0)&0\\
0&C
\end{bmatrix} \cdot
\begin{bmatrix} (\delta w,\delta a)\\
\delta b
\end{bmatrix},
$$
where $C$ is a $\ssc$-isomorphism. An $\ssc$-operator of this form
is $\ssc$-Fredholm if and only if the linearization $f'(0,0)$ is
sc-Fredholm. In that case the Fredholm indices are the same. This
completes the proof.
\end{proof}

 What we discussed in this paper is a minimal set of
concepts needed to develop a Fredholm theory. The next paper will
contain a treatment of implicit function theorems in the splicing
context.

We refer the reader to the upcoming volume \cite{HWZ-polyfolds1} for
a more exhaustive list of splicing constructions and constructs.

\end{document}